\documentclass{article}%
\usepackage{amsmath}
\usepackage{amsfonts}
\usepackage{amssymb}
\usepackage{graphicx}%
\providecommand{\U}[1]{\protect \rule{.1in}{.1in}}
\newtheorem{theorem}{Theorem}

\newtheorem{definition}[theorem]{Definition}
\newtheorem{example}[theorem]{Example}

\newtheorem{lemma}[theorem]{Lemma}

\newtheorem{proposition}[theorem]{Proposition}
\newtheorem{remark}[theorem]{Remark}

\newenvironment{proof}[1][Proof]{\noindent \textbf{#1.} }{\  \rule{0.5em}{0.5em}}
\begin{document}

\title{Multi-Dimensional $G$--Brownian Motion and Related Stochastic Calculus under
$G$--Expectation}
\author{Shige PENG\thanks{The author thanks the partial support from the Natural
Science Foundation of China, grant No. 10131040. He thanks to the anonymous
referee of SPA for his constructive suggestions. In particular, the
discussions on Jensen's inequlity for $G$--convex functions in Section 6 were
stimulated by one of his interesting questions. }\\Institute of Mathematics\\Shandong University\\250100, Jinan, China\\peng@sdu.edu.cn}
\date{1st version: arXiv:math.PR/0601699 v1 28 Jan 2006 \\
This version: November 9, 2006}
\maketitle

\noindent \textbf{Abstract. }{\small We develop a notion of nonlinear
expectation ----}$G${\small --expectation---- generated by a nonlinear heat
equation with infinitesimal generator }$G${\small . We first study
multi-dimensional }$G${\small --normal distributions. With this nonlinear
distribution we can introduce our }$G${\small --expectation under which the
canonical process is a multi--dimensional }$G${\small --Brownian motion. We
then establish the related stochastic calculus, especially stochastic
integrals of It\^{o}'s type with respect to our }$G${\small --Brownian motion
and derive the related It\^{o}'s formula. We have also obtained the existence
and uniqueness of stochastic differential equation under our }$G$%
{\small --expectation. }\newline


\noindent \textbf{MSC 2000 Classification Numbers: }60H10, 60H05, 60H30, 35K55,
35K15, 49L25\newline \newline

\section{Introduction}

The purpose of this paper is to extend classical stochastic calculus for
multi-dimensional Brownian motion to the setting of nonlinear $G$%
--expectation. We first recall the general framework of nonlinear expectation
studied in \cite{Peng2005} and \cite{Peng2004}, where the usual linearity is
replaced by positive homogeneity and subadditivity. Such a sublinear
expectation functional enables us to construct a Banach space, similar to an
$\mathbb{L}^{1}$-space, starting from a functional lattice of Daniell's type.

Then we proceed to construct a sublinear expectation on the space of
continuous paths from $\mathbb{R}_{+}$ to $\mathbb{R}^{d}$, starting from $0$,
which will be an analogue of Wiener's law. The operation mainly consists in
replacing the Brownian semigroup by a nonlinear semigroup coming from the
solution of a nonlinear parabolic partial differential equation (\ref{eq-heat}%
) where appears a mapping $G$ acting on Hessian matrices. Indeed, the Markov
property permits to define in the same way nonlinear conditional expectations
with respect to the past. Then we presents some rules and examples of
computations under the newly constructed $G$-Brownian (motion) expectation.
The fact that the underlying marginal nonlinear expectations are $G$-normal
distributions derived from the nonlinear heat equation (\ref{eq-heat}) is very
helpful to estimate natural functionals. As result, our $G$--Brownian motion
also has independent increments with identical $G$-normal distributions.

$G$--Brownian motion has a very rich and interesting new structure which non
trivially generalizes the classical one. We thus can establish the related
stochastic calculus, especially $G$--It\^{o}'s integrals (see \cite[1942]%
{Ito}) and the related quadratic variation process $\left \langle
B\right \rangle $. A very interesting new phenomenon of our $G$-Brownian motion
is that its quadratic process $\left \langle B\right \rangle $ also has
independent increments which are identically distributed. \ The corresponding
$G$--It\^{o}'s formula is obtained. We then introduce the notion of
$G$--martingales and the related Jensen inequality for a new type of
\textquotedblleft$G$--convex\textquotedblright \ functions. We have also
established the existence and uniqueness of solution to stochastic
differential equation under our stochastic calculus by the same Picard
iterations as in the classical situation. Books on stochastic calculus e.g.,
\cite{Chu-Will}, \cite{HWY},\  \cite{IW}, \cite{Ito-McKean}, \cite{KSh},
\cite{Oksendal}, \cite{Protter}, \cite{Revuz-Yor}, \cite{Yong-Zhou} are
recommended for understanding the present results and some further possible
developments of this new stochastic calculus.

As indicated in Remark \ref{Rem-1}, the nonlinear expectations discussed in
this paper can be regarded as coherent risk measures. This with the related
conditional expectations $\mathbb{E}[\cdot|\mathcal{H}_{t}]_{t\geq0}$ makes a
dynamic risk measure: $G$--risk measure.

The other motivation of our $G$--expectation is the notion of (nonlinear)
$g$--expectations introduced in \cite{Peng1997}, \cite{Peng1997b}. Here $g$ is
the generating function of a backward stochastic differential equation (BSDE)
on a given probability space $(\Omega,\mathcal{F},\mathbf{P})$. The natural
definition of the conditional $g$--expectations with respect to the past
induces rich properties of nonlinear $g$--martingale theory (see among others,
\cite{BCHMP1}, \cite{Chen98}, \cite{CE}, \cite{CKJ}, \cite{CHMP},
\cite{CHMP3}, \cite{CP}, \cite{CP1}, \cite{Jiang}, \cite{JC}, \cite{Peng1999},
\cite{Peng2005a}, \cite{Peng2005b}, \cite{PX2003}). Recently $g$--expectations
are also studied as dynamic risk measures: $g$--risk measure (cf.
\cite{Roazza2003}, \cite{El-Bar}, \cite{DPR}). Fully nonlinear super-hedging
is also a possible application (cf. \cite{Lyons} and \cite{Touzi} where new
BSDE\ approach was introduced).

The notion of $g$--expectation is defined on a given probability space. In
\cite{Peng2005} (see also \cite{Peng2004}), we have constructed a kind of
filtration--consistent nonlinear expectations through the so--called nonlinear
Markov chain. As compared with the framework of $g$--expectations, the theory
of $G$--expectation is intrinsic, a meaning similar to \textquotedblleft
intrinsic geometry\textquotedblright \ in the sense that it is not based on a
given (linear) probability space. \ Since the classical Brownian expectation
as well as many other linear and nonlinear expectations are dominated by our
$G$--Expectation (see Remark \ref{Rem-17}, Example \ref{Exa-AB} and
\cite{Peng2005}) and thus can be considered as continuous functionals, our
theory also provides a flexible theoretical framework.

$1$--dimensional $G$--Brownian motion was studied in \cite{Peng2006a}. Unlike
the classical situation, in general, we cannot find a system of coordinates
under which the corresponding components $B^{i}$, $i=1,\cdots,d$, are mutually
independent from each others. The mutual quadratic variations $\left \langle
B^{i},B^{j}\right \rangle $ will play essentially important rule.

During the reversion process of this paper, the author has found a very
interesting paper \cite{Denis-M} by Denis and Martini on super-pricing of
contingent claims under model uncertainty of volatility. They have introduced
a norm on the space of continuous paths $\Omega=C([0,T])$ which corresponds to
our $L_{G}^{2}$-norm and developed a stochastic integral. There is no notions
of nonlinear expectation such as $G$--expectation, conditional $G$%
--expectation, the related $G$-normal distribution and the notion of
independence in their paper. But on the other hand, powerful tools in capacity
theory enables them to obtain pathwise results of random variables and
stochastic processes through the language of \textquotedblleft
quasi-surely\textquotedblright, (see Feyel and de La Pradelle \cite{Feyel-DLP}%
) in the place of \textquotedblleft almost surely\textquotedblright \ in
classical probability theory. Their method provides a way to proceed a
pathwise analysis for our $G$--Brownian motion and the related stochastic
calculus under $G$--expectation, see our forthcoming paper joint with Denis.

This paper is organized as follows: in Section 2, we recall the framework of
nonlinear expectation established in \cite{Peng2005} and adapt it to our
objective. In section 3 we introduce $d$--dimensional $G$-normal distribution
and discuss its main properties. In Section 4 we introduce $d$--dimensional
$G$-Brownian motion, the corresponding $G$--expectation and their main
properties. We then can establish stochastic integral with respect to
$G$-Brownian motion of It\^{o}'s type, the related quadratic variation
processes and then $G$--It\^{o}'s formula in Section 5, $G$--martingale and
the Jensen's inequality for $G$--convex functions in Section 6 and the
existence and uniqueness theorem of SDE driven by $G$-Brownian motion in
Section 7.

The whole results of this paper are based on the very basic knowledge of
Banach space and the parabolic partial differential equation (\ref{eq-heat}).
When this $G$-heat equation (\ref{eq-heat}) is linear, our $G$-Brownian motion
becomes the classical Brownian motion. This paper still provides an analytic
shortcut to reach the sophistic It\^{o}'s calculus.

\section{Nonlinear expectation: a general framework}

We briefly recall the notion of nonlinear expectations introduced in
\cite{Peng2005}. Following Daniell's famous integration (cf. Daniell 1918
\cite{Daniell}, see also \cite{Yan}), we begin with a vector lattice. Let
$\Omega$ be a given set and let $\mathcal{H}$ be a vector lattice of real
functions defined on $\Omega$ containing $1$, namely, $\mathcal{H}$ is a
linear space such that $1\in \mathcal{H}$ and that $X\in \mathcal{H}$ implies
$|X|\in \mathcal{H}$. $\mathcal{H}$ is a space of random variables. We assume
the functions on $\mathcal{H}$ are all bounded.

\begin{definition}
\label{Def-1}\textbf{A nonlinear expectation }$\mathbb{E}$ is a functional
$\mathcal{H}\mapsto \mathbb{R}$ satisfying the following properties\newline%
\newline \textbf{(a) Monotonicity:} if $X,Y\in \mathcal{H}$ and $X\geq Y$ then
$\mathbb{E}[X]\geq \mathbb{E}[Y].$\newline \textbf{(b)} \textbf{Preserving of
constants:} $\mathbb{E}[c]=c$.\newline \newline In this paper we are interested
in the sublinear expectations which satisfy\newline \newline \textbf{(c)}
\textbf{Sub-additivity (or self--dominated property):}%
\[
\mathbb{E}[X]-\mathbb{E}[Y]\leq \mathbb{E}[X-Y],\  \  \forall X,Y\in \mathcal{H}.
\]
\textbf{(d) Positive homogeneity: } $\mathbb{E}[\lambda X]=\lambda
\mathbb{E}[X]$,$\  \  \forall \lambda \geq0$, $X\in \mathcal{H}$.\newline%
\textbf{(e) Constant translatability: }$\mathbb{E}[X+c]=\mathbb{E}[X]+c$.
\end{definition}

\medskip

\begin{remark}
It is clear that (d)+(e) implies (b). \label{Rem-1}We recall the notion of the
above sublinear expectations was systematically introduced by Artzner,
Delbaen, Eber and Heath \cite{ADEH1}, \cite{ADEH2}, in the case where $\Omega$
is a finite set, and by Delbaen \cite{Delbaen} in general situation with the
notation of risk measure: $\rho(X)=\mathbb{E}[-X]$. See also in Huber
\cite{Huber} for even early study of this notion $\mathbb{E}$ (called upper
expectation $\mathbf{E}^{\ast}$ in Ch.10 of \cite{Huber}).
\end{remark}

We follow \cite{Peng2005} to introduce a Banach space via $\mathcal{H}$ and
$\mathbb{E}$. We denote $\left \Vert X\right \Vert :=\mathbb{E}[|X|]$,
$X\in \mathcal{H}$. $\mathcal{H}$ forms a normed space $(\mathcal{H},\left \Vert
\cdot \right \Vert )$ under $\left \Vert \cdot \right \Vert $ in the following
sense. For each $X,Y\in \mathcal{H}$ such that $\left \Vert X-Y\right \Vert =0$,
we set $X=Y$. This is equivalent to say that the linear subspace
\[
\mathcal{H}_{0}:=\{X\in \mathcal{H},\  \left \Vert X\right \Vert =0\}
\]
is the null space, or in other words, we only consider the elements in the
quotient space $\mathcal{H}/\mathcal{H}_{0}$. Under such arrangement
$(\mathcal{H},\left \Vert \cdot \right \Vert )$ is a normed space. We denote by
$([\mathcal{H}],\left \Vert \cdot \right \Vert )$, or simply $[\mathcal{H}]$, the
completion of $(\mathcal{H},\left \Vert \cdot \right \Vert )$. $(\mathcal{H}%
,\left \Vert \cdot \right \Vert )$ is a dense subspace of the Banach space
$([\mathcal{H}],\left \Vert \cdot \right \Vert )$ (see e.g., Yosida \cite{Yosida}
Sec. I-10).

For any $X\in \mathcal{H}$, the mappings
\begin{align*}
X^{+}(\omega)  &  =\max \{X(\omega),0\}:\mathcal{H\longmapsto H},\\
X^{-}(\omega)  &  =\max \{-X(\omega),0\}:\mathcal{H\longmapsto H}%
\end{align*}
satisfy
\begin{align*}
|X^{+}-Y^{+}|  &  \leq|X-Y|,\\
X^{-}-Y^{-}  &  \leq(Y-X)^{+}\leq|X-Y|.\
\end{align*}
Thus they are both contract mappings under $\left \Vert \cdot \right \Vert $ and
can be continuously extended to the Banach space $[\mathcal{H}]$.

We define the partial order \textquotedblleft$\geq$\textquotedblright \ in this
Banach space.

\begin{definition}
An element $X$ in $([\mathcal{H}],\left \Vert \cdot \right \Vert )$ is said to be
nonnegative, or $X\geq0$, $0\leq X$, if $X=X^{+}$. We also denote by $X\geq
Y$, or $Y\leq X$. if $X-Y\geq0$.
\end{definition}

It is easy to check that if $X\geq Y$ and $Y\geq X$, then $X=Y$ in
$([\mathcal{H}],\left \Vert \cdot \right \Vert )$. The nonlinear expectation
$\mathbb{E}[\cdot]$ can be continuously extended to $([\mathcal{H}],\left \Vert
\cdot \right \Vert )$ on which \textbf{(a)--(e)} still hold.

\section{$G$--normal distributions}

For a given positive integer $n$, we will denote by $(x,y)$ the scalar product
of $x$, $y\in \mathbb{R}^{n}$ and by $\left \vert x\right \vert =(x,x)^{1/2}$ the
Euclidean norm of $x$. We denote by $lip(\mathbb{R}^{n})$ the space of all
bounded and Lipschitz real functions on $\mathbb{R}^{n}$. We introduce the
notion of nonlinear distribution-- $G$--normal distribution. A $G$%
\textbf{--normal distribution} is a nonlinear expectation defined on
$lip(\mathbb{R}^{d})$ (here $\mathbb{R}^{d}$ is considered as $\Omega$ and
$lip(\mathbb{R}^{d})$ as $\mathcal{H})$:
\[
P_{1}^{G}(\phi)=u(1,0):\phi \in lip(\mathbb{R}^{d})\mapsto \mathbb{R}%
\]
where $u=u(t,x)$ is a bounded continuous function on $[0,\infty)\times
\mathbb{R}^{d}$ which is the viscosity solution of the following nonlinear
parabolic partial differential equation (PDE)
\begin{equation}
\frac{\partial u}{\partial t}-G(D^{2}u)=0,\  \ u(0,x)=\phi(x),\ (t,x)\in
\lbrack0,\infty)\times \mathbb{R}^{d}, \label{eq-heat}%
\end{equation}
here $D^{2}u$ is the Hessian matrix of $u$, i.e., $D^{2}u=(\partial
_{x^{i}x^{j}}^{2}u)_{i,j=1}^{d}$ and
\begin{equation}
G(A)=G_{\Gamma}(A)=\frac{1}{2}\sup_{\gamma \in \Gamma}\text{tr}[\gamma \gamma
^{T}A],\  \ A=(A_{ij})_{i,j=1}^{d}\in \mathbb{S}_{d}. \label{GammaA}%
\end{equation}
$\mathbb{S}_{d}$ denotes the space of $d\times d$ symmetric matrices. $\Gamma$
is a given non empty, bounded and closed subset of $\mathbb{R}^{d\times d}$,
the space of all $d\times d$ matrices.

\begin{remark}
\label{Rem-3}The nonlinear heat equation (\ref{eq-heat}) is a special kind of
Hamilton--Jacobi--Bellman equation. The existence and uniqueness of
(\ref{eq-heat}) in the sense of viscosity solution can be found in, for
example, \cite{CIL}, \cite{FS}, \cite{Peng1992}, \cite{Yong-Zhou}, and
\cite{Krylov} for $C^{1,2}$-solution if $\gamma \gamma^{T}\geq \sigma_{0}I_{n}$,
for each $\gamma \in \Gamma$, for a given constant $\sigma_{0}>0$ (see also in
\cite{Oksendal} for elliptic cases). It is a known result that $u(t,\cdot)\in
lip(\mathbb{R}^{d})$ (see e.g. \cite{Yong-Zhou} Ch.4, prop.3.1. or
\cite{Peng1992} Lemma 3.1. for the Lipschitz continuity of $u(t,\cdot)$, or
Lemma 5.5 and Proposition 5.6 in \cite{Peng2004} for a more general
conclusion). The boundedness is simply from the comparison theorem (or maximum
principle) of this PDE. It is also easy to check that, for a given $\psi \in
lip(\mathbb{R}^{d}\times \mathbb{R}^{d})$, $P_{1}^{G}(\psi(x,\cdot))$ is still
a bounded and Lipschitz function in $x$.
\end{remark}

In the case where $\Gamma$ is a singleton $\{ \gamma_{0}\}$, the above PDE
becomes a standard linear heat equation and thus, for $G^{0}=G_{\{ \gamma
_{0}\}}$, the corresponding $G^{0}$--distribution is just the $d$--dimensional
classical normal distribution $\mathcal{N}(0,\gamma_{0}\gamma_{0}^{T})$. In a
typical case where $\gamma_{0}=I_{d}\in \Gamma$, we have
\[
P_{1}^{G^{0}}(\phi)=\frac{1}{(2\pi)^{d/2}}\int_{\mathbb{R}^{d}}\exp
[-\sum_{i=1}^{d}\frac{(x^{i})^{2}}{2}]\phi(x)dx.
\]

In the case where $\gamma_{0}\in \Gamma$, from comparison theorem of PDE,
\begin{equation}
P_{1}^{G}(\phi)\geq P_{1}^{G^{0}}(\phi),\  \forall \phi \in lip(\mathbb{R}%
^{d})\text{.}\  \label{compar}%
\end{equation}
More generally, for each subset $\Gamma^{\prime}\subset \Gamma$, the
corresponding $P^{G_{\Gamma^{\prime}}}$--distribution is dominated by $P^{G}$
in the following sense:%
\[
P_{1}^{G_{\Gamma^{\prime}}}(\phi)-P_{1}^{G_{\Gamma^{\prime}}}(\psi)\leq
P_{1}^{G}(\phi-\psi),\  \  \forall \phi,\psi \in lip(\mathbb{R}^{d}).
\]

\begin{remark}
In \cite{Peng2006a} we have discussed $1$--dimensional case, which corresponds
$d=1$ and $\Gamma=[\sigma,1]\subset \mathbb{R}$, where $\sigma \in \lbrack0,1]$
is a given constant. In this case the nonlinear heat equation (\ref{eq-heat})
becomes%
\[
\frac{\partial u}{\partial t}-\frac{1}{2}[(\partial_{xx}^{2}u)^{+}-\sigma
^{2}(\partial_{xx}^{2}u)^{-}]=0,\  \ u(0,x)=\phi(x),\ (t,x)\in \lbrack
0,\infty)\times \mathbb{R}.
\]
In multi--dimensional case we also have the following typical nonlinear heat
equation:%
\[
\frac{\partial u}{\partial t}-\frac{1}{2}\sum_{i=1}^{d}[(\partial_{x^{i}x^{i}%
}^{2}u)^{+}-\sigma_{i}^{2}(\partial_{x^{i}x^{i}}^{2}u)^{-}]=0
\]
where $\sigma_{i}\in \lbrack0,1]$ are given constants. In this case we have
\[
\Gamma=\{diag[\gamma_{1},\cdots,\gamma_{d}],\gamma_{i}\in \lbrack \sigma
_{i},1],\ i=1,\cdots,d\}.
\]

\end{remark}

The corresponding normal distribution with mean at $x\in \mathbb{R}^{d}$ and
square variation $t>0$ is $P_{1}^{G}(\phi(x+\sqrt{t}\times \cdot))$. Just like
the classical situation of a normal distribution,, we have

\begin{lemma}
\label{Scaling}For each $\phi \in lip(\mathbb{R}^{d})$, the function
\begin{equation}
u(t,x)=P_{1}^{G}(\phi(x+\sqrt{t}\times \cdot)),\  \ (t,x)\in \lbrack
0,\infty)\times \mathbb{R}^{d} \label{u(t,x)}%
\end{equation}
is the solution of the nonlinear heat equation (\ref{eq-heat}) with the
initial condition $u(0,\cdot)=\phi(\cdot)$.
\end{lemma}

\begin{proof}
Let $u\in C([0,\infty)\times \mathbb{R}^{d})$ be the viscosity solution of
(\ref{eq-heat}) with $u(0,\cdot)=\phi(\cdot)\in lip(\mathbb{R}^{d})$. For a
fixed $(\bar{t},\bar{x})\in(0,\infty)\times \mathbb{R}^{d}$, we denote $\bar
{u}(t,x)=u(t\times \bar{t},x\sqrt{\bar{t}}+\bar{x})$. Then $\bar{u}$ is the
viscosity solution of (\ref{eq-heat}) with the initial condition $\bar
{u}(0,x)=\phi(x\sqrt{\bar{t}}+\bar{x})$. Indeed, let $\psi$ be a $C^{1,2}$
function on $(0,\infty)\times \mathbb{R}^{d}$ such that $\psi \geq \bar{u}$
(resp. $\psi \leq \bar{u}$) and $\psi(\tau,\xi)=\bar{u}(\tau,\xi)$ for a fixed
$(\tau,\xi)\in(0,\infty)\times \mathbb{R}^{d}$. We have $\psi(\frac{t}{\bar{t}%
},\frac{x-\bar{x}}{\sqrt{\bar{t}}})\geq u(t,x)$, for all $(t,x)$ and
\[
\psi(\frac{t}{\bar{t}},\frac{x-\bar{x}}{\sqrt{\bar{t}}})=u(t,x),\  \text{at
}(t,x)=(\tau \bar{t},\xi \sqrt{\bar{t}}+\bar{x}).
\]
Since $u$ is the viscosity solution of (\ref{eq-heat}), at the point
$(t,x)=(\tau \bar{t},\xi \sqrt{\bar{t}}+\bar{x})$, we have
\[
\frac{\partial \psi(\frac{t}{\bar{t}},\frac{x-\bar{x}}{\sqrt{\bar{t}}}%
)}{\partial t}-G(D^{2}\psi(\frac{t}{\bar{t}},\frac{x-\bar{x}}{\sqrt{\bar{t}}%
}))\leq0\  \ (\text{resp. }\geq0).
\]
But $G$ is a positive homogenous function, i.e., $G(\lambda A)=\lambda G(A)$,
when $\lambda \geq0$, we thus derive
\[
\frac{\partial \psi(t,x)}{\partial t}-G(D^{2}\psi(t,x))|_{(t,x)=(\tau,\xi)}%
\leq0\  \ (\text{resp. }\geq0).
\]
This implies that $\bar{u}$ is the viscosity subsolution (resp. supersolution)
of (\ref{eq-heat}). According to the definition of $P^{G}(\cdot)$ we obtain
(\ref{u(t,x)}).
\end{proof}

\begin{definition}
\label{Def-2}We denote
\begin{equation}
P_{t}^{G}(\phi)(x)=P_{1}^{G}(\phi(x+\sqrt{t}\times \cdot))=u(t,x),\  \ (t,x)\in
\lbrack0,\infty)\times \mathbb{R}^{d}. \label{P_tx}%
\end{equation}

\end{definition}

From the above lemma, for each $\phi \in lip(\mathbb{R}^{d})$, we have the
following nonlinear version of chain rule:%
\begin{equation}
P_{t}^{G}(P_{s}^{G}(\phi))(x)=P_{t+s}^{G}(\phi)(x),\  \ s,t\in \lbrack
0,\infty),\ x\in \mathbb{R}^{d}.\  \label{Chapman}%
\end{equation}
This chain rule was initialled by Nisio \cite{Nisio1} and \cite{Nisio2} in
terms of \textquotedblleft envelope of Markovian semi-groups\textquotedblright%
. See also \cite{Peng2005}.

\begin{lemma}
\label{a-x}The solution of (\ref{eq-heat}) with initial condition
$u(0,x)=\phi((\mathbf{a},x))$, for a given $\phi \in lip(\mathbb{R})$, has the
form $u(t,x)=\bar{u}(t,\bar{x})$, $\bar{x}=(\mathbf{a},x)$, where $\bar{u}$ is
the solution of
\begin{equation}
\frac{\partial \bar{u}}{\partial t}-G_{\mathbf{a}}(\partial_{\bar{x}\bar{x}%
}\bar{u})=0,\  \ u(0,\bar{x})=\phi(\bar{x}),\ (t,\bar{x})\in \lbrack
0,\infty)\times \mathbb{R}, \label{PDE-G}%
\end{equation}
where
\[
G_{\mathbf{a}}(\beta)=\frac{1}{2}\max_{\gamma \in \Gamma}tr[\gamma \gamma
^{T}\mathbf{aa}^{T}\beta],\  \  \beta \in \mathbb{R}.
\]
The above PDE can be written
\begin{equation}
\frac{\partial \bar{u}}{\partial t}-\frac{1}{2}[\sigma_{\mathbf{aa}^{T}%
}(\partial_{\bar{x}\bar{x}}\bar{u})^{+}+\sigma_{-\mathbf{aa}^{T}}%
(\partial_{\bar{x}\bar{x}}\bar{u})^{-}]=0,\  \ u(0,\bar{x})=\phi(\bar{x}).
\label{PDE-Gg}%
\end{equation}
where we denote $\mathbf{aa}^{T}=[a^{i}a^{j}]_{i,j=1}^{d}\in \mathbb{S}_{d}$
and
\begin{equation}
\sigma_{A}=\sup_{\gamma \in \Gamma}tr[\gamma \gamma^{T}A]=2G(A),\  \  \ A\in
\mathbb{S}_{d}. \label{sigma}%
\end{equation}
Here $\mathbb{S}_{d}$ is the space of $d\times d$ symmetric matrices.
\end{lemma}

\begin{remark}
\label{Rem-a-x}It is clear that the functional
\[
P_{1}^{G_{\mathbf{a}}}(\phi)=\bar{u}(1,0):\phi \in lip(\mathbb{R}%
)\mapsto \mathbb{R}%
\]
constitutes a special $1$--dimensional nonlinear normal distribution, called
$G_{\mathbf{a}}$--\textbf{normal distribution.}
\end{remark}

\begin{proof}
It is clear that the PDE (\ref{PDE-G}) has a unique viscosity solution. We
then can set $u(t,x)=\bar{u}(t,(\mathbf{a},x))$ and check that $u$ is the
viscosity solution of (\ref{eq-heat}). (\ref{PDE-Gg}) is then easy to check.
\end{proof}

\begin{example}
\label{Convexfunc}In the above lemma, if $\phi$ is convex, and $\sigma
_{\mathbf{aa}^{T}}>0$, then
\[
P_{t}^{G}(\phi((\mathbf{a},\cdot))(x)=\frac{1}{\sqrt{2\pi \sigma_{\mathbf{aa}%
^{T}}t}}\int_{-\infty}^{\infty}\phi(y)\exp(-\frac{(y-x)^{2}}{2\sigma
_{\mathbf{aa}^{T}}t})dy.
\]
If $\phi$ is concave and $\sigma_{-\mathbf{aa}^{T}}<0$, then%
\[
P_{t}^{G}(\phi((\mathbf{a},\cdot))(x)=\frac{1}{\sqrt{2\pi|\sigma
_{-\mathbf{aa}^{T}}|t}}\int_{-\infty}^{\infty}\phi(y)\exp(-\frac{(y-x)^{2}%
}{2|\sigma_{-\mathbf{aa}^{T}}|t})dy.
\]

\end{example}

\begin{proposition}
\label{Prop-BM9}We have \newline(i) For each $t>0$, the $G$--normal
distribution $P_{t}^{G}$ is a nonlinear expectation on the lattice
$lip(\mathbb{R}^{d})$, with $\Omega=\mathbb{R}^{d}$, satisfying (a)--(e) of
definition \ref{Def-1}. The corresponding \ completion space $[\mathcal{H]=[}%
lip(\mathbb{R}^{d})]_{t}$ under the norm $\left \Vert \phi \right \Vert
_{t}:=P_{t}^{G}(|\phi|)(0)$ contains $\phi(x)=x_{1}^{n_{1}}\times \cdots \times
x_{d}^{n_{d}}$, $n_{i}=1,2,\cdots$, $i=1,\cdots,d$, $x=(x_{1},\cdots
,x_{d})^{T}$ as well as $x_{1}^{n_{1}}\times \cdots \times x_{d}^{n_{d}}%
\times \psi(x)$, $\psi \in lip(\mathbb{R}^{d}\mathbb{)}$ as its special
elements. Relation (\ref{P_tx}) still holds. We also have the following
properties\newline(ii) We have, for each $\mathbf{a}=(a_{1},\cdots,a_{d}%
)^{T}\in \mathbb{R}^{d}$ and $A\in \mathbb{S}_{d}$
\begin{align*}
P_{t}^{G}((\mathbf{a},x)_{x\in \mathbb{R}^{d}})  &  =0,\  \  \  \  \  \  \\
P_{t}^{G}(((\mathbf{a},x)^{2})_{x\in \mathbb{R}^{d}})  &  =t\cdot
\sigma_{\mathbf{aa}^{T}},\  \ P_{t}^{G}((-(\mathbf{a},x)^{2})_{x\in
\mathbb{R}^{d}})=t\cdot \sigma_{-\mathbf{aa}^{T}},\  \\
P_{t}^{G}(((\mathbf{a},x)^{4})_{x\in \mathbb{R}^{d}})  &  =6(\sigma
_{\mathbf{aa}^{T}})t^{2},\  \ P_{t}^{G}((-(\mathbf{a},x)^{4})_{x\in
\mathbb{R}^{d}})=-6(\sigma_{-\mathbf{aa}^{T}})^{2}t^{2},\\
P_{t}^{G}(((Ax,x))_{x\in \mathbb{R}^{d}})  &  =t\cdot \sigma_{A}=2G(A)t.
\end{align*}

\end{proposition}

\begin{proof}
(ii) By Lemma \ref{a-x}, we have the explicit solutions of the nonlinear PDE
(\ref{eq-heat}) with the following different initial condition $u(0,x)=\phi
(x)$:
\[%
\begin{array}
[c]{ccl}%
\phi(x)=(\mathbf{a},x) & \Longrightarrow & u(t,x)=(\mathbf{a},x),\\
\phi(x)=(\mathbf{a},x)^{4} & \Longrightarrow & u(t,x)=(\mathbf{a}%
,x)^{4}+6(\mathbf{a},x)^{2}\sigma_{\mathbf{aa}^{T}}t+6\sigma_{\mathbf{aa}^{T}%
}^{2}t^{2},\\
\phi(x)=-(\mathbf{a},x)^{4} & \Longrightarrow & u(t,x)=-(\mathbf{a}%
,x)^{4}+6(\mathbf{a},x)^{2}\sigma_{-\mathbf{aa}^{T}}t-6|\sigma_{-\mathbf{aa}%
^{T}}|^{2}t^{2}.
\end{array}
\]
Similarly, we can check that $%
\begin{array}
[c]{ccl}%
\phi(x)=(Ax,x) & \Longrightarrow & u(t,x)=(Ax,x)+\sigma_{A}t\text{. }%
\end{array}
$This implies, by setting $A=\mathbf{aa}^{T}$ and $A=-\mathbf{aa}^{T}$,
\[%
\begin{array}
[c]{ccl}%
\phi(x)=(\mathbf{a},x)^{2} & \Longrightarrow & u(t,x)=(\mathbf{a}%
,x)^{2}+\sigma_{\mathbf{aa}^{T}}t,\\
\phi(x)=-(\mathbf{a},x)^{2} & \Longrightarrow & u(t,x)=-(\mathbf{a}%
,x)^{2}+\sigma_{-\mathbf{aa}^{T}}t.
\end{array}
\]
More generally, for $\phi(x)=(\mathbf{a},x)^{2n}$, we have
\[
u(t,x)=\frac{1}{\sqrt{2\pi \sigma_{\mathbf{aa}^{T}}t}}\int_{-\infty}^{\infty
}y^{2n}\exp(-\frac{(y-x)^{2}}{2\sigma_{\mathbf{aa}^{T}}t})dy.
\]
By this we can prove (i).
\end{proof}

\section{$G$--Brownian motions under $G$--expectations}

In the rest of this paper, we set $\Omega=C_{0}^{d}(\mathbb{R}^{+})$ the space
of all $\mathbb{R}^{d}$--valued continuous paths $(\omega_{t})_{t\in
\mathbb{R}^{+}}$, with $\omega_{0}=0$, equipped with the distance
\[
\rho(\omega^{1},\omega^{2}):=\sum_{i=1}^{\infty}2^{-i}[(\max_{t\in \lbrack
0,i]}|\omega_{t}^{1}-\omega_{t}^{2}|)\wedge1].
\]
$\mathbf{\Omega}$ is the classical canonical space and $\omega=(\omega
_{t})_{t\geq0}$ is the corresponding canonical process. It is well--known that
in this canonical space there exists a Wiener measure $(\mathbf{\Omega
},\mathcal{F},P)$ under which the canonical process $B_{t}(\omega)=\omega_{t}$
is a $d$--dimensional Brownian motion.

For each fixed $T\geq0$, we consider the following space of random variables:
\[
L_{ip}^{0}(\mathcal{H}_{T}):=\{X(\omega)=\phi(\omega_{t_{1}},\cdots
,\omega_{t_{m}}),\forall m\geq1,\;t_{1},\cdots,t_{m}\in \lbrack0,T],\phi \in
l_{ip}(\mathbb{R}^{d\times m})\}.
\]
It is clear that $\{L_{ip}^{0}(\mathcal{H}_{t})\}_{t\geq0}$ constitute a
family of sub-lattices such that $L_{ip}^{0}(\mathcal{H}_{t})\subseteq
L_{ip}^{0}(\mathcal{H}_{T})$, for $t\leq T<\infty$. $L_{ip}^{0}(\mathcal{H}%
_{t})$ representing the past history of $\omega$ at the time $t$. It's
completion will play the same role of Brownian filtration $\mathcal{F}_{t}%
^{B}$ as in classical stochastic analysis. We also denote%
\[
L_{ip}^{0}(\mathcal{H}):=%
{\displaystyle \bigcup \limits_{n=1}^{\infty}}
L_{ip}^{0}(\mathcal{H}_{n}).
\]

\begin{remark}
It is clear that $lip(\mathbb{R}^{d\times m})$ and then $L_{ip}^{0}%
(\mathcal{H}_{T})$, $L_{ip}^{0}(\mathcal{H})$ are vector lattices. Moreover,
since $\phi,\psi \in lip(\mathbb{R}^{d\times m})$ implies $\phi \cdot \psi \in
lip(\mathbb{R}^{d\times m})$ thus $X$, $Y\in L_{ip}^{0}(\mathcal{H}_{T})$
implies $X\cdot Y\in L_{ip}^{0}(\mathcal{H}_{T})$; $X$, $Y\in L_{ip}%
^{0}(\mathcal{H})$ implies $X\cdot Y\in L_{ip}^{0}(\mathcal{H})$.
\end{remark}

We will consider the canonical space and set $B_{t}(\omega)=\omega_{t}$,
$t\in \lbrack0,\infty)$, for $\omega \in \Omega$.

\begin{definition}
\label{Def-3}The canonical process $B$ is called a ($d$--dimensional)
$G$\textbf{--Brownian} \textbf{motion} under a nonlinear expectation
$\mathbb{E}$ defined on $L_{ip}^{0}(\mathcal{H})$ if \newline \textrm{(i) }For
each $s,t\geq0$ and $\psi \in lip(\mathbb{R}^{d})$, $B_{t}$ and $B_{t+s}-B_{s}$
are identically distributed:
\[
\mathbb{E}[\psi(B_{t+s}-B_{s})]=\mathbb{E}[\psi(B_{t})]=P_{t}^{G}(\psi).\  \
\]
\textrm{(ii) }For each $m=1,2,\cdots$, $0\leq t_{1}<\cdots<t_{m}<\infty$, the
increment $B_{t_{m}}-B_{t_{m-1}}$ is \textquotedblleft
backwardly\textquotedblright \ independent from $B_{t_{1}}$,$\cdots,B_{t_{m-1}%
}$ in the following sense: for each $\phi \in lip(\mathbb{R}^{d\times m})$,%
\[
\mathbb{E}[\phi(B_{t_{1}},\cdots,B_{t_{m-1}},B_{t_{m}})]=\mathbb{E}[\phi
_{1}(B_{t_{1}},\cdots,B_{t_{m-1}})],
\]
where $\phi_{1}(x^{1},\cdots,x^{m-1})=\mathbb{E}[\phi(x^{1},\cdots
,x^{m-1},B_{t_{m}}-B_{t_{m-1}}+x^{m-1})]$, $x^{1}$,$\cdots,x^{m-1}%
\in \mathbb{R}^{d}$. \newline The related conditional expectation of
$\phi(B_{t_{1}},\cdots,B_{t_{m}})$ under $\mathcal{H}_{t_{k}}$ is defined by%
\begin{equation}
\mathbb{E}[\phi(B_{t_{1}},\cdots,B_{t_{k}},\cdots,B_{t_{m}})|\mathcal{H}%
_{t_{k}}]=\phi_{m-k}(B_{t_{1}},\cdots,B_{t_{k}}), \label{Condition}%
\end{equation}
where
\[
\phi_{m-k}(x^{1},\cdots,x^{k})=\mathbb{E}[\phi(x^{1},\cdots,x^{k},B_{t_{k+1}%
}-B_{t_{k}}+x^{k},\cdots,B_{t_{m}}-B_{t_{k}}+x^{k})].
\]

\end{definition}

It is proved in \cite{Peng2005} that $\mathbb{E}[\cdot]$ consistently defines
a nonlinear expectation on the vector lattice $L_{ip}^{0}(\mathcal{H}_{T})$ as
well as on $L_{ip}^{0}(\mathcal{H})$ satisfying (a)--(e) in Definition
\ref{Def-1}. It follows that $\mathbb{E}[|X|]$, $X\in L_{ip}^{0}%
(\mathcal{H}_{T})$ (resp. $L_{ip}^{0}(\mathcal{H})$) forms a norm and thus
$L_{ip}^{0}(\mathcal{H}_{T})$ (resp. $L_{ip}^{0}(\mathcal{H})$) can be
extended, under this norm, to a Banach space. We denote this space by
$L_{G}^{1}(\mathcal{H}_{T})$ (resp. $L_{G}^{1}(\mathcal{H})$). For each $0\leq
t\leq T<\infty$, we have $L_{G}^{1}(\mathcal{H}_{t})\subseteq L_{G}%
^{1}(\mathcal{H}_{T})\subset L_{G}^{1}(\mathcal{H})$. In $L_{G}^{1}%
(\mathcal{H}_{T})$ (resp. $L_{G}^{1}(\mathcal{H}_{T})$), $\mathbb{E}[\cdot]$
still satisfies (a)--(e) in Definition \ref{Def-1}.

\begin{remark}
It is suggestive to denote $L_{ip}^{0}(\mathcal{H}_{t})$ by $\mathcal{H}%
_{t}^{0}$ and $L_{G}^{1}(\mathcal{H}_{t})$ by $\mathcal{H}_{t}$, $L_{G}%
^{1}(\mathcal{H})$ by $\mathcal{H}$ and thus consider the conditional
expectation $\mathbb{E}[\cdot|\mathcal{H}_{t}]$ as a projective mapping from
$\mathcal{H}$ to $\mathcal{H}_{t}$. The notation $L_{G}^{1}(\mathcal{H}_{t})$
is due to the similarity of $L^{1}(\Omega,\mathcal{F}_{t},P)$ in classical
stochastic analysis.
\end{remark}

\begin{definition}
The expectation $\mathbb{E}[\cdot]:L_{G}^{1}(\mathcal{H})\mapsto \mathbb{R}$
introduced through above procedure is called $G$\textbf{--expectation}, or
$G$--Brownian expectation. The corresponding canonical process $B$ is said to
be a $G$--Brownian motion under $\mathbb{E}[\cdot]$.
\end{definition}

For a given $p>1$, we also denote $L_{G}^{p}(\mathcal{H})=\{X\in L_{G}%
^{1}(\mathcal{H}),\ |X|^{p}\in L_{G}^{1}(\mathcal{H})\}$. $L_{G}%
^{p}(\mathcal{H})$ is also a Banach space under the norm $\left \Vert
X\right \Vert _{p}:=(\mathbb{E}[|X|^{p}])^{1/p}$. We have (see Appendix)
\[
\left \Vert X+Y\right \Vert _{p}\leq \left \Vert X\right \Vert _{p}+\left \Vert
Y\right \Vert _{p}%
\]
and, for each $X\in L_{G}^{p}$, $Y\in L_{G}^{q}(Q)$ with $\frac{1}{p}+\frac
{1}{q}=1$,%
\[
\left \Vert XY\right \Vert =\mathbb{E}[|XY|]\leq \left \Vert X\right \Vert
_{p}\left \Vert X\right \Vert _{q}.
\]
With this we have $\left \Vert X\right \Vert _{p}\leq \left \Vert X\right \Vert
_{p^{\prime}}$ if $p\leq p^{\prime}$.

We now consider the conditional expectation introduced in (\ref{Condition}).
For each fixed $t=t_{k}\leq T$, the conditional expectation $\mathbb{E}%
[\cdot|\mathcal{H}_{t}]:L_{ip}^{0}(\mathcal{H}_{T})\mapsto L_{ip}%
^{0}(\mathcal{H}_{t})$ is a continuous mapping under $\left \Vert
\cdot \right \Vert $. Indeed, we have $\mathbb{E}[\mathbb{E}[X|\mathcal{H}%
_{t}]]=\mathbb{E}[X]$, $X\in L_{ip}^{0}(\mathcal{H}_{T})$ and, since
$P_{t}^{G}$ is subadditive,
\[
\mathbb{E}[X|\mathcal{H}_{t}]-\mathbb{E}[Y|\mathcal{H}_{t}]\leq \mathbb{E}%
[X-Y|\mathcal{H}_{t}]\leq \mathbb{E}[|X-Y||\mathcal{H}_{t}]
\]
We thus obtain%
\[
\mathbb{E[E}[X|\mathcal{H}_{t}]-\mathbb{E}[Y|\mathcal{H}_{t}]]\leq
\mathbb{E}[X-Y]
\]
and
\[
\left \Vert \mathbb{E}[X|\mathcal{H}_{t}]-\mathbb{E}[Y|\mathcal{H}%
_{t}]\right \Vert \leq \left \Vert X-Y\right \Vert .
\]
It follows that $\mathbb{E}[\cdot|\mathcal{H}_{t}]$ can be also extended as a
continuous mapping $L_{G}^{1}(\mathcal{H}_{T})\mapsto L_{G}^{1}(\mathcal{H}%
_{t})$. If the above $T$ is not fixed, then we can obtain $\mathbb{E}%
[\cdot|\mathcal{H}_{t}]:L_{G}^{1}(\mathcal{H})\mapsto L_{G}^{1}(\mathcal{H}%
_{t})$.

\begin{proposition}
\label{Prop-1-7}We list the properties of $\mathbb{E}[\cdot|\mathcal{H}_{t}]$,
$t\in \lbrack0,T]$, that hold in $L_{ip}^{0}(\mathcal{H}_{T})$ and still hold
for $X$, $Y\in$ $L_{G}^{1}(\mathcal{H}_{T})$:\newline \newline \textbf{(i)
}$\mathbb{E}[X|\mathcal{H}_{t}]=X$, for $X\in L_{G}^{1}(\mathcal{H}_{t})$,
$t\leq T$.\newline \textbf{(ii) }If $X\geq Y$, then $\mathbb{E}[X|\mathcal{H}%
_{t}]\geq \mathbb{E}[Y|\mathcal{H}_{t}]$.\newline \textbf{(iii) }$\mathbb{E}%
[X|\mathcal{H}_{t}]-\mathbb{E}[Y|\mathcal{H}_{t}]\leq \mathbb{E}%
[X-Y|\mathcal{H}_{t}].$\newline \textbf{(iv) }$\mathbb{E}[\mathbb{E}%
[X|\mathcal{H}_{t}]|\mathcal{H}_{s}]=\mathbb{E}[X|\mathcal{H}_{t\wedge s}]$,
$\mathbb{E}[\mathbb{E}[X|\mathcal{H}_{t}]]=\mathbb{E}[X].$\newline \textbf{(v)}
$\mathbb{E}[X+\eta|\mathcal{H}_{t}]=\mathbb{E}[X|\mathcal{H}_{t}]+\eta$,
$\eta \in L_{G}^{1}(\mathcal{H}_{t})$\newline \textbf{(vi)} $\mathbb{E}[\eta
X|\mathcal{H}_{t}]=\eta^{+}\mathbb{E}[X|\mathcal{H}_{t}]+\eta^{-}%
\mathbb{E}[-X|\mathcal{H}_{t}]$, for bounded $\eta \in L_{G}^{1}(\mathcal{H}%
_{t}).$\newline \textbf{(vii)} We have the following independence:
\[
\mathbb{E}[X|\mathcal{H}_{t}]=\mathbb{E}[X],\medskip \  \  \forall X\in L_{G}%
^{1}(\mathcal{H}_{T}^{t}),\  \forall T\geq0,
\]
\  \newline where $L_{G}^{1}(\mathcal{H}_{T}^{t})$ is the extension, under
$\left \Vert \cdot \right \Vert $, of $L_{ip}^{0}(\mathcal{H}_{T}^{t})$ which
consists of random variables of the form $\phi(B_{t_{1}}^{t},B_{t_{2}}%
^{t},\cdots,B_{t_{m}}^{t})$, $\phi \in lip(\mathbb{R}^{m})$, $t_{1}%
,\cdots,t_{m}\in \lbrack0,T]$, $m=1,2,\cdots$. Here we denote%
\[
B_{s}^{t}=B_{t+s}-B_{t},\  \ s\geq0.
\]
\newline \textbf{(viii)} The increments of $B$ are identically distributed:%
\[
\mathbb{E}[\phi(B_{t_{1}}^{t},B_{t_{2}}^{t},\cdots,B_{t_{m}}^{t}%
)]=\mathbb{E}[\phi(B_{t_{1}},B_{t_{2}},\cdots,B_{t_{m}})].
\]

\end{proposition}

The meaning of the independence in (vii) is similar to the classical one:

\begin{definition}
An $\mathbb{R}^{n}$ valued random variable $Y\in(L_{G}^{1}(\mathcal{H}))^{n}$
is said to be independent of $\mathcal{H}_{t}$ for some given $t$ if for each
$\phi \in lip(\mathbb{R}^{n})$ we have%
\[
\mathbb{E}[\phi(Y)|\mathcal{H}_{t}]=\mathbb{E}[\phi(Y)].
\]

\end{definition}

It is seen that the above property (vii) also holds for the situation $X\in$
$L_{G}^{1}(\mathcal{H}^{t})$ where $L_{G}^{1}(\mathcal{H}^{t})$ is the
completion of the sub-lattice $\cup_{T\geq0}L_{G}^{1}(\mathcal{H}_{T}^{t})$
under $\left \Vert \cdot \right \Vert $.

From the above results we have

\begin{proposition}
For each fixed $t\geq0$, $(B_{s}^{t})_{s\geq0}$ is a $G$--Brownian motion in
$L_{G}^{1}(\mathcal{H}^{t})$ under the same $G$--expectation $\mathbb{E}%
[\cdot]$.
\end{proposition}

\begin{remark}
We can also prove, using Lemma \ref{Scaling}, that the time scaling of $B$,
i.e., $\tilde{B}=(\sqrt{\lambda}B_{t/\lambda})_{t\geq0}$ also consists a
$G$--Brownian motion.
\end{remark}

The following property is very useful

\begin{proposition}
\label{E-x+y}Let $X,Y\in L_{G}^{1}(\mathcal{H})$ be such that $\mathbb{E}%
[Y|\mathcal{H}_{t}]=-\mathbb{E}[-Y|\mathcal{H}_{t}]$, for some $t\in
\lbrack0,T]$. Then we have%
\[
\mathbb{E}[X+Y|\mathcal{H}_{t}]=\mathbb{E}[X|\mathcal{H}_{t}]+\mathbb{E}%
[Y|\mathcal{H}_{t}].
\]
In particular, if $\mathbb{E}[Y|\mathcal{H}_{t}]=\mathbb{E}[-Y|\mathcal{H}%
_{t}]=0$, then $\mathbb{E}[X+Y|\mathcal{H}_{t}]=\mathbb{E}[X|\mathcal{H}_{t}]$.
\end{proposition}

\begin{proof}
It is simply because we have $\mathbb{E}[X+Y|\mathcal{H}_{t}]\leq
\mathbb{E}[X|\mathcal{H}_{t}]+\mathbb{E}[Y|\mathcal{H}_{t}]$ and
\[
\mathbb{E}[X+Y|\mathcal{H}_{t}]\geq \mathbb{E}[X|\mathcal{H}_{t}]-\mathbb{E}%
[-Y|\mathcal{H}_{t}]=\mathbb{E}[X|\mathcal{H}_{t}]+\mathbb{E}[Y|\mathcal{H}%
_{t}]\text{.}%
\]

\end{proof}

\begin{example}
\label{Exm-GBM-14a}From the last relation of Proposition \ref{Prop-BM9}-(ii),
we have%
\[
\mathbb{E}[(AB_{t},B_{t})]=\sigma_{A}t=2G(A)t,\  \  \forall A\in \mathbb{S}_{d}.
\]
More general, for each $s\leq t$ and $\eta=(\eta^{ij})_{i,j=1}^{d}\in
L_{G}^{2}(\mathcal{H}_{s};\mathbb{S}_{d})$,
\begin{equation}
\mathbb{E}[(\eta B_{t}^{s},B_{t}^{s})|\mathcal{H}_{s}]=\sigma_{\eta}%
t=2G(\eta)t,\ s,t\geq0. \label{eq-GMB-14a}%
\end{equation}

\end{example}

\begin{definition}
We will denote, in the rest of this paper,
\begin{equation}
B_{t}^{\mathbf{a}}= ( \mathbf{a},B_{t}) ,\  \  \  \text{for each }\mathbf{a}%
=(a_{1},\cdots,a_{d})^{T}\in \mathbb{R}^{d} \label{a-Brown}%
\end{equation}
From Lemma \ref{a-x} and Remark \ref{Rem-a-x},
\[
\mathbb{E}[\phi(B_{t}^{\mathbf{a}})]=P_{t}^{G}(\phi( ( \mathbf{a},\cdot)
))=P_{t}^{G_{\mathbf{a}}}(\phi)
\]
where $P^{G_{\mathbf{a}}}$ is the ($1$--dimensional) $G_{\mathbf{a}}$--normal
distribution. Thus, according to Definition \ref{Def-3} for $d$--dimensional
$G$--Brownian motion, $B^{\mathbf{a}}$ forms a $1$--dimensional $G_{\mathbf{a}%
}$--\textbf{Brownian motion} for which the $G_{\mathbf{a}}$--expectation
coincides with $\mathbb{E}[\cdot]$.
\end{definition}

\begin{example}
\label{Exam-1}For each $0\leq s-t$, we have
\[
\mathbb{E}[\psi(B_{t}-B_{s})|\mathcal{H}_{s}]=\mathbb{E}[\psi(B_{t}-B_{s})]
\]
If $\phi$ is a real convex function on $\mathbb{R}$ and at least not growing
too fast, then%
\begin{align*}
&  \mathbb{E}[X\phi(B_{T}^{\mathbf{a}}-B_{t}^{\mathbf{a}})|\mathcal{H}_{t}]\\
&  =X^{+}\mathbb{E}[\phi(B_{T}^{\mathbf{a}}-B_{t}^{\mathbf{a}})|\mathcal{H}%
_{t}]+X^{-}\mathbb{E}[-\phi(B_{T}^{\mathbf{a}}-B_{t}^{\mathbf{a}}%
)|\mathcal{H}_{t}]\\
&  =\frac{X^{+}}{\sqrt{2\pi(T-t)\sigma_{\mathbf{aa}^{T}}}}\int_{-\infty
}^{\infty}\phi(x)\exp(-\frac{x^{2}}{2(T-t)\sigma_{\mathbf{aa}^{T}}})dx\\
&  -\frac{X^{-}}{\sqrt{2\pi(T-t)|\sigma_{-\mathbf{aa}^{T}}|}}\int_{-\infty
}^{\infty}\phi(x)\exp(-\frac{x^{2}}{2(T-t)|\sigma_{-\mathbf{aa}^{T}}|})dx.
\end{align*}
In particular, for $n=1,2,\cdots,$%
\begin{align*}
\mathbb{E}[|B_{t}^{\mathbf{a}}-B_{s}^{\mathbf{a}}|^{n}|\mathcal{H}_{s}]  &
=\mathbb{E}[|B_{t-s}^{\mathbf{a}}|^{n}]\\
&  =\frac{1}{\sqrt{2\pi(t-s)\sigma_{\mathbf{aa}^{T}}}}\int_{-\infty}^{\infty
}|x|^{n}\exp(-\frac{x^{2}}{2(t-s)\sigma_{\mathbf{aa}^{T}}})dx.
\end{align*}
But we have $\mathbb{E}[-|B_{t}^{\mathbf{a}}-B_{s}^{\mathbf{a}}|^{n}%
|\mathcal{H}_{s}]=\mathbb{E}[-|B_{t-s}^{\mathbf{a}}|^{n}]$ which is $0$ when
$\sigma_{-\mathbf{aa}^{T}}=0$ and
\[
\frac{-1}{\sqrt{2\pi(t-s)|\sigma_{-\mathbf{aa}^{T}}|}}\int_{-\infty}^{\infty
}|x|^{n}\exp(-\frac{x^{2}}{2(t-s)|\sigma_{-\mathbf{aa}^{T}}|})dx,\  \  \text{if
}\sigma_{-\mathbf{aa}^{T}}<0.\
\]
Exactly as in classical cases, we have $\mathbb{E}[B_{t}^{\mathbf{a}}%
-B_{s}^{\mathbf{a}}|\mathcal{H}_{s}]=0$ and
\begin{align*}
\mathbb{E}[(B_{t}^{\mathbf{a}}-B_{s}^{\mathbf{a}})^{2}|\mathcal{H}_{s}]  &
=\sigma_{\mathbf{aa}^{T}}(t-s),\  \  \  \mathbb{E}[(B_{t}^{\mathbf{a}}%
-B_{s}^{\mathbf{a}})^{4}|\mathcal{H}_{s}]=3\sigma_{\mathbf{aa}^{T}}%
^{2}(t-s)^{2},\\
\mathbb{E}[(B_{t}^{\mathbf{a}}-B_{s}^{\mathbf{a}})^{6}|\mathcal{H}_{s}]  &
=15\sigma_{\mathbf{aa}^{T}}^{3}(t-s)^{3},\  \  \mathbb{E}[(B_{t}^{\mathbf{a}%
}-B_{s}^{\mathbf{a}})^{8}|\mathcal{H}_{s}]=105\sigma_{\mathbf{aa}^{T}}%
^{4}(t-s)^{4},\\
\mathbb{E}[|B_{t}^{\mathbf{a}}-B_{s}^{\mathbf{a}}||\mathcal{H}_{s}]  &
=\frac{\sqrt{2(t-s)\sigma_{\mathbf{aa}^{T}}}}{\sqrt{\pi}},\  \  \mathbb{E}%
[|B_{t}^{\mathbf{a}}-B_{s}^{\mathbf{a}}|^{3}|\mathcal{H}_{s}]=\frac{2\sqrt
{2}[(t-s)\sigma_{\mathbf{aa}^{T}}]^{3/2}}{\sqrt{\pi}},\\
\mathbb{E}[|B_{t}^{\mathbf{a}}-B_{s}^{\mathbf{a}}|^{5}|\mathcal{H}_{s}]  &
=8\frac{\sqrt{2}[(t-s)\sigma_{\mathbf{aa}^{T}}]^{5/2}}{\sqrt{\pi}}.
\end{align*}

\end{example}

\begin{example}
\label{Exam-2}For each $n=1,2,\cdots,$ $0\leq t\leq T$ and $X\in L_{G}%
^{1}(\mathcal{H}_{t})$, we have%
\[
\mathbb{E}[X(B_{T}^{\mathbf{a}}-B_{t}^{\mathbf{a}})|\mathcal{H}_{t}%
]=X^{+}\mathbb{E}[(B_{T}^{\mathbf{a}}-B_{t}^{\mathbf{a}})|\mathcal{H}%
_{t}]+X^{-}\mathbb{E}[-(B_{T}^{\mathbf{a}}-B_{t}^{\mathbf{a}})|\mathcal{H}%
_{t}]=0.
\]
This with Proposition \ref{E-x+y} yields%
\[
\mathbb{E}[Y+X(B_{T}^{\mathbf{a}}-B_{t}^{\mathbf{a}})|\mathcal{H}%
_{t}]=\mathbb{E}[Y|\mathcal{H}_{t}],\  \ Y\in L_{G}^{1}(\mathcal{H}).
\]
We also have,
\begin{align*}
\mathbb{E}[X(B_{T}^{\mathbf{a}}-B_{t}^{\mathbf{a}})^{2}|\mathcal{H}_{t}]  &
=X^{+}\mathbb{E}[(B_{T}^{\mathbf{a}}-B_{t}^{\mathbf{a}})^{2}|\mathcal{H}%
_{t}]+X^{-}\mathbb{E}[-(B_{T}^{\mathbf{a}}-B_{t}^{\mathbf{a}})^{2}%
|\mathcal{H}_{t}]\\
&  =[X^{+}\sigma_{\mathbf{aa}^{T}}+X^{-}\sigma_{-\mathbf{aa}^{T}}](T-t).
\end{align*}

\end{example}

\begin{remark}
\label{Rem-17}It is clear that we can define an expectation $E[\cdot]$ on
$L_{ip}^{0}(\mathcal{H})$ in the same way as in Definition \ref{Def-3} with
the standard normal distribution $P_{1}^{0}(\cdot)$ in the place of $P_{1}%
^{G}(\cdot)$. If $I_{d}\in \Gamma$, then it follows from (\ref{compar}) that
$P_{1}^{0}(\cdot)$ is dominated by $P_{1}^{G}(\cdot)$ in the sense
\[
P_{1}^{0}(\phi)-P_{1}^{0}(\psi)\leq P_{1}^{G}(\phi-\psi).
\]
Then $E[\cdot]$ can be continuously extended to $L_{G}^{1}(\mathcal{H})$.
$E[\cdot]$ is a linear expectation under which $(B_{t})_{t\geq0}$ behaves as a
Brownian motion. We have
\begin{equation}
-\mathbb{E}[-X]\leq E^{0}[X]\leq \mathbb{E}[X],\  \ -\mathbb{E}[-X|\mathcal{H}%
_{t}]\leq E^{0}[X|\mathcal{H}_{t}]\leq \mathbb{E}[X|\mathcal{H}_{t}].
\end{equation}
More generally, if $\Gamma^{\prime}\subset \Gamma$, since the corresponding
$P^{\prime}=P^{G_{\Gamma^{\prime}}}$ is dominated by $P^{G}=P^{G_{\Gamma}}$,
thus the corresponding expectation $\mathbb{E}^{\prime}$ is well--defined in
$L_{G}^{1}(\mathcal{H})$ and $\mathbb{E}^{\prime}$ is dominated by
$\mathbb{E}$:
\[
\mathbb{E}^{\prime}[X]-\mathbb{E}^{\prime}[Y]\leq \mathbb{E}[X-Y],\  \ X,Y\in
L_{G}^{1}(\mathcal{H}).
\]
Such kind of extension through the above type of domination relations was
discussed in details in \cite{Peng2005}. With this domination we then can
introduce a large kind of time consistent linear or nonlinear expectations and
the corresponding conditional expectations, not necessarily to be positive
homogeneous and/or subadditive, as continuous functionals in $L_{G}%
^{1}(\mathcal{H})$. See Example \ref{Exa-AB} for a further discussion.
\end{remark}

\begin{example}
\label{Exam-B2}Since
\[
\mathbb{E}[2B_{s}^{\mathbf{a}}(B_{t}^{\mathbf{a}}-B_{s}^{\mathbf{a}%
})|\mathcal{H}_{s}]=\mathbb{E}[-2B_{s}^{\mathbf{a}}(B_{t}^{\mathbf{a}}%
-B_{s}^{\mathbf{a}})|\mathcal{H}_{s}]=0,
\]
we have,
\begin{align*}
\mathbb{E}[(B_{t}^{\mathbf{a}})^{2}-(B_{s}^{\mathbf{a}})^{2}|\mathcal{H}_{s}]
&  =\mathbb{E}[(B_{t}^{\mathbf{a}}-B_{s}^{\mathbf{a}}+B_{s}^{\mathbf{a}}%
)^{2}-(B_{s}^{\mathbf{a}})^{2}|\mathcal{H}_{s}]\\
&  =\mathbb{E}[(B_{t}^{\mathbf{a}}-B_{s}^{\mathbf{a}})^{2}+2(B_{t}%
^{\mathbf{a}}-B_{s}^{\mathbf{a}})B_{s}^{\mathbf{a}}|\mathcal{H}_{s}]\\
&  =\sigma_{\mathbf{aa}^{T}}(t-s)
\end{align*}
and%
\begin{align*}
\mathbb{E}[((B_{t}^{\mathbf{a}})^{2}-(B_{s}^{\mathbf{a}})^{2})^{2}%
|\mathcal{H}_{s}]  &  =\mathbb{E}[\{(B_{t}^{\mathbf{a}}-B_{s}^{\mathbf{a}%
}+B_{s}^{\mathbf{a}})^{2}-(B_{s}^{\mathbf{a}})^{2}\}^{2}|\mathcal{H}_{s}]\\
&  =\mathbb{E}[\{(B_{t}^{\mathbf{a}}-B_{s}^{\mathbf{a}})^{2}+2(B_{t}%
^{\mathbf{a}}-B_{s}^{\mathbf{a}})B_{s}^{\mathbf{a}}\}^{2}|\mathcal{H}_{s}]\\
&  =\mathbb{E}[(B_{t}^{\mathbf{a}}-B_{s}^{\mathbf{a}})^{4}+4(B_{t}%
^{\mathbf{a}}-B_{s}^{\mathbf{a}})^{3}B_{s}^{\mathbf{a}}+4(B_{t}^{\mathbf{a}%
}-B_{s}^{\mathbf{a}})^{2}(B_{s}^{\mathbf{a}})^{2}|\mathcal{H}_{s}]\\
&  \leq \mathbb{E}[(B_{t}^{\mathbf{a}}-B_{s}^{\mathbf{a}})^{4}]+4\mathbb{E}%
[|B_{t}^{\mathbf{a}}-B_{s}^{\mathbf{a}}|^{3}]|B_{s}^{\mathbf{a}}%
|+4\sigma_{\mathbf{aa}^{T}}(t-s)(B_{s}^{\mathbf{a}})^{2}\\
&  =3\sigma_{\mathbf{aa}^{T}}^{2}(t-s)^{2}+8\sqrt{\frac{2}{\pi}}%
[\sigma_{\mathbf{aa}^{T}}(t-s)]^{3/2}|B_{s}^{\mathbf{a}}|+4\sigma
_{\mathbf{aa}^{T}}(t-s)(B_{s}^{\mathbf{a}})^{2}%
\end{align*}

\end{example}

\section{It\^{o}'s integral of $G$--Brownian motion}

\subsection{Bochner's integral}

\begin{definition}
\label{Def-4}For $T\in \mathbb{R}_{+}$, a partition $\pi_{T}$ of $[0,T]$ is a
finite ordered subset $\pi=\{t_{1},\cdots,t_{N}\}$ such that $0=t_{0}%
<t_{1}<\cdots<t_{N}=T$.
\[
\mu(\pi_{T})=\max \{|t_{i+1}-t_{i}|,i=0,1,\cdots,N-1\} \text{.}%
\]
We use $\pi_{T}^{N}=\{t_{0}^{N}<t_{1}^{N}<\cdots<t_{N}^{N}\}$ to denote a
sequence of partitions of $[0,T]$ such that $\lim_{N\rightarrow \infty}\mu
(\pi_{T}^{N})=0$.
\end{definition}

Let $p\geq1$ be fixed. We consider the following type of simple processes: for
a given partition $\{t_{0},\cdots,t_{N}\}=\pi_{T}$ of $[0,T]$, we set%
\[
\eta_{t}(\omega)=\sum_{k=0}^{N-1}\xi_{k}(\omega)\mathbf{I}_{[t_{k},t_{k+1}%
)}(t)
\]
where $\xi_{k}\in L_{G}^{p}(\mathcal{H}_{t_{i}})$, $k=0,1,2,\cdots,N-1$ are
given. The collection of these type of processes is denoted by $M_{G}%
^{p,0}(0,T)$.

\begin{definition}
\label{Def-5}For an $\eta \in M_{G}^{1,0}(0,T)$ with $\eta_{t}=\sum_{k=0}%
^{N-1}\xi_{k}(\omega)\mathbf{I}_{[t_{k},t_{k+1})}(t)$, the related Bochner
integral is
\[
\int_{0}^{T}\eta_{t}(\omega)dt=\sum_{k=0}^{N-1}\xi_{k}(\omega)(t_{k+1}%
-t_{k}).
\]

\end{definition}

\begin{remark}
We set, for each $\eta \in M_{G}^{1,0}(0,T)$,
\[
\mathbb{\tilde{E}}_{T}[\eta]:=\frac{1}{T}\int_{0}^{T}\mathbb{E}[\eta
_{t}]dt=\frac{1}{T}\sum_{k=0}^{N-1}\mathbb{E}\xi_{k}(\omega)(t_{k+1}-t_{k}).
\]
It is easy to check that $\mathbb{\tilde{E}}_{T}:M_{G}^{1,0}(0,T)\longmapsto
\mathbb{R}$ forms a nonlinear expectation satisfying (a)--(e) of Definition
\ref{Def-1}. We then can introduce a nature norm
\[
\left \Vert \eta \right \Vert _{T}^{1}=\mathbb{\tilde{E}}_{T}[|\eta|]=\frac{1}%
{T}\int_{0}^{T}\mathbb{E}[|\eta_{t}|]dt.
\]
Under this norm $M_{G}^{1,0}(0,T)$ can extended to $M_{G}^{1}(0,T)$ which is a
Banach space.
\end{remark}

\begin{definition}
For each $p\geq1$, we denote by $M_{G}^{p}(0,T)$ the completion of
$M_{G}^{p,0}(0,T)$ under the norm%
\[
(\frac{1}{T}\int_{0}^{T}\left \Vert |\eta_{t}|^{p}\right \Vert dt)^{1/p}=\left(
\frac{1}{T}\sum_{k=0}^{N-1}\mathbb{E[}|\xi_{k}(\omega)|^{p}](t_{k+1}%
-t_{k})\right)  ^{1/p}.
\]

\end{definition}

We observe that,
\begin{equation}
\mathbb{E}[|\int_{0}^{T}\eta_{t}(\omega)dt|]\leq \sum_{k=0}^{N-1}\left \Vert
\xi_{k}(\omega)\right \Vert (t_{k+1}-t_{k})=\int_{0}^{T}\mathbb{E}[|\eta
_{t}|]dt. \label{Bohner}%
\end{equation}
We then have

\begin{proposition}
The linear mapping $\int_{0}^{T}\eta_{t}(\omega)dt:M_{G}^{1,0}(0,T)\mapsto
L_{G}^{1}(\mathcal{H}_{T})$ is continuous and thus can be continuously
extended to $M_{G}^{1}(0,T)\mapsto L_{G}^{1}(\mathcal{H}_{T})$. We still
denote this extended mapping by $\int_{0}^{T}\eta_{t}(\omega)dt$, $\eta \in
M_{G}^{1}(0,T)$.
\end{proposition}

Since $M_{G}^{p}(0,T)\subset M_{G}^{1}(0,T)$, for $p\geq1$. Thus this
definition holds for $\eta \in M_{G}^{p}(0,T)$.

\subsection{It\^{o}'s integral of $G$--Brownian motion}

We still use $B_{t}^{\mathbf{a}}:= (\mathbf{a},B_{t})$ as in (\ref{a-Brown}).

\begin{definition}
For each $\eta \in M_{G}^{2,0}(0,T)$ with the form $\eta_{t}(\omega)=\sum
_{k=0}^{N-1}\xi_{k}(\omega)\mathbf{I}_{[t_{k},t_{k+1})}(t)$, we define
\[
I(\eta)=\int_{0}^{T}\eta(s)dB_{s}^{\mathbf{a}}:=\sum_{k=0}^{N-1}\xi
_{k}(B_{t_{k+1}}^{\mathbf{a}}-B_{t_{k}}^{\mathbf{a}})\mathbf{.}%
\]

\end{definition}

\begin{lemma}
\label{bdd}The linear mapping $I:M_{G}^{2,0}(0,T)\longmapsto L_{G}%
^{2}(\mathcal{H}_{T})$ is continuous and thus can be continuously extended to
$I:M_{G}^{2}(0,T)\longmapsto L_{G}^{2}(\mathcal{H}_{T})$. In fact we have, for
each $\eta \in M_{G}^{2,0}(0,T)$,
\begin{align}
\mathbb{E}[\int_{0}^{T}\eta(s)dB_{s}^{\mathbf{a}}]  &  =0,\  \  \label{e1}\\
\mathbb{E}[(\int_{0}^{T}\eta(s)dB_{s}^{\mathbf{a}})^{2}]  &  \leq
\sigma_{\mathbf{aa}^{T}}\int_{0}^{T}\mathbb{E}[\eta^{2}(s)]ds. \label{e2}%
\end{align}

\end{lemma}

\begin{definition}
We define, for a fixed $\eta \in M_{G}^{2}(0,T)$ the stochastic calculus
\[
\int_{0}^{T}\eta(s)dB_{s}^{\mathbf{a}}:=I(\eta).
\]
It is clear that (\ref{e1}), (\ref{e2}) still hold for $\eta \in M_{G}%
^{2}(0,T)$.
\end{definition}

\textbf{Proof of Lemma \ref{bdd}. }From Example \ref{Exam-2}, for each $k$,
\[
\mathbb{E}\mathbf{[}\xi_{k}(B_{t_{k+1}}^{\mathbf{a}}-B_{t_{k}}^{\mathbf{a}%
})|\mathcal{H}_{t_{k}}]=0.
\]
We have%
\begin{align*}
\mathbb{E}[\int_{0}^{T}\eta(s)dB_{s}^{\mathbf{a}}]  &  =\mathbb{E[}\int
_{0}^{t_{N-1}}\eta(s)dB_{s}^{\mathbf{a}}+\xi_{N-1}(B_{t_{N}}^{\mathbf{a}%
}-B_{t_{N-1}}^{\mathbf{a}})]\\
&  =\mathbb{E[}\int_{0}^{t_{N-1}}\eta(s)dB_{s}^{\mathbf{a}}+\mathbb{E}%
\mathbf{[}\xi_{N-1}(B_{t_{N}}^{\mathbf{a}}-B_{t_{N-1}}^{\mathbf{a}%
})|\mathcal{H}_{t_{N-1}}]]\\
&  =\mathbb{E[}\int_{0}^{t_{N-1}}\eta(s)dB_{s}^{\mathbf{a}}].
\end{align*}
We then can repeat this procedure to obtain (\ref{e1}). We now prove
(\ref{e2})
\begin{align*}
&  \mathbb{E}[\left(  \int_{0}^{T}\eta(s)dB_{s}^{\mathbf{a}}\right)
^{2}]=\mathbb{E[}\left(  \int_{0}^{t_{N-1}}\eta(s)dB_{s}^{\mathbf{a}}%
+\xi_{N-1}(B_{t_{N}}^{\mathbf{a}}-B_{t_{N-1}}^{\mathbf{a}})\right)  ^{2}]\\
&  =\mathbb{E[}\left(  \int_{0}^{t_{N-1}}\eta(s)dB_{s}^{\mathbf{a}}\right)
^{2}\\
&  +\mathbb{E}[2\left(  \int_{0}^{t_{N-1}}\eta(s)dB_{s}^{\mathbf{a}}\right)
\xi_{N-1}(B_{t_{N}}^{\mathbf{a}}-B_{t_{N-1}}^{\mathbf{a}})+\xi_{N-1}%
^{2}(B_{t_{N}}^{\mathbf{a}}-B_{t_{N-1}}^{\mathbf{a}})^{2}|\mathcal{H}%
_{t_{N-1}}]]\\
&  =\mathbb{E[}\left(  \int_{0}^{t_{N-1}}\eta(s)dB_{s}^{\mathbf{a}}\right)
^{2}+\xi_{N-1}^{2}\sigma_{\mathbf{aa}^{T}}(t_{N}-t_{N-1})].
\end{align*}
Thus $\mathbb{E}[(\int_{0}^{t_{N}}\eta(s)dB_{s}^{\mathbf{a}})^{2}%
]\leq \mathbb{E[}\left(  \int_{0}^{t_{N-1}}\eta(s)dB_{s}^{\mathbf{a}}\right)
^{2}]+\mathbb{E}[\xi_{N-1}^{2}]\sigma_{\mathbf{aa}^{T}}(t_{N}-t_{N-1})$. We
then repeat this procedure to deduce
\[
\mathbb{E}[(\int_{0}^{T}\eta(s)dB_{s})^{2}]\leq \sigma_{\mathbf{aa}^{T}}%
\sum_{k=0}^{N-1}\mathbb{E}[(\xi_{k})^{2}](t_{k+1}-t_{k})=\int_{0}%
^{T}\mathbb{E}[(\eta(t))^{2}]dt.
\]
$\blacksquare$

We list some main property of the It\^{o}'s integral of $G$--Brownian motion.
We denote for some $0\leq s\leq t\leq T$,
\[
\int_{s}^{t}\eta_{u}dB_{u}^{\mathbf{a}}:=\int_{0}^{T}\mathbf{I}_{[s,t]}%
(u)\eta_{u}dB_{u}^{\mathbf{a}}.
\]
We have

\begin{proposition}
\label{Prop-Integ}Let $\eta,\theta \in M_{G}^{2}(0,T)$ and let $0\leq s\leq
r\leq t\leq T$. Then in $L_{G}^{1}(\mathcal{H}_{T})$ we have\newline(i)
$\int_{s}^{t}\eta_{u}dB_{u}^{\mathbf{a}}=\int_{s}^{r}\eta_{u}dB_{u}%
^{\mathbf{a}}+\int_{r}^{t}\eta_{u}dB_{u}^{\mathbf{a}}.$\newline(ii) $\int
_{s}^{t}(\alpha \eta_{u}+\theta_{u})dB_{u}^{\mathbf{a}}=\alpha \int_{s}^{t}%
\eta_{u}dB_{u}^{\mathbf{a}}+\int_{s}^{t}\theta_{u}dB_{u}^{\mathbf{a}}%
,\ $if$\  \alpha$ is bounded and in $L_{G}^{1}(\mathcal{H}_{s})$,\newline(iii)
$\mathbb{E[}X+\int_{r}^{T}\eta_{u}dB_{u}^{\mathbf{a}}|\mathcal{H}%
_{s}]=\mathbb{E[}X|\mathcal{H}_{s}]$, $\forall X\in L_{G}^{1}(\mathcal{H}%
)$\newline(iv) $\mathbb{E[}(\int_{r}^{T}\eta_{u}dB_{u}^{\mathbf{a}}%
)^{2}|\mathcal{H}_{s}]\leq \sigma_{\mathbf{aa}^{T}}\int_{r}^{T}\mathbb{E[}%
|\eta_{u}|^{2}|\mathcal{H}_{s}]du$
\end{proposition}

\subsection{Quadratic variation process of $G$--Brownian motion}

We now consider the quadratic variation of $G$--Brownian motion. It
concentrically reflects the characteristic of the `uncertainty' part of the
$G$-Brownian motion $B$. This makes a major difference from the classical
Brownian motion.

Let $\pi_{t}^{N}$, $N=1,2,\cdots$, be a sequence of partitions of $[0,t]$. We consider%

\begin{align*}
(B_{t}^{\mathbf{a}})^{2}  &  =\sum_{k=0}^{N-1}[(B_{t_{k+1}^{N}}^{\mathbf{a}%
})^{2}-(B_{t_{k}^{N}}^{\mathbf{a}})^{2}]\\
&  =\sum_{k=0}^{N-1}2B_{t_{k}^{N}}^{\mathbf{a}}(B_{t_{k+1}^{N}}^{\mathbf{a}%
}-B_{t_{k}^{N}}^{\mathbf{a}})+\sum_{k=0}^{N-1}(B_{t_{k+1}^{N}}^{\mathbf{a}%
}-B_{t_{k}^{N}}^{\mathbf{a}})^{2}%
\end{align*}
As $\mu(\pi_{t}^{N})=\max_{0\leq k\leq N-1}(t_{k+1}^{N}-t_{k}^{N}%
)\rightarrow0$, the first term of the right side tends to $\int_{0}^{t}%
B_{s}^{\mathbf{a}}dB_{s}^{\mathbf{a}}$. The second term must converge. We
denote its limit by $\left \langle B^{\mathbf{a}}\right \rangle _{t}$, i.e.,
\begin{equation}
\left \langle B^{\mathbf{a}}\right \rangle _{t}=\lim_{\mu(\pi_{t}^{N}%
)\rightarrow0}\sum_{k=0}^{N-1}(B_{t_{k+1}^{N}}^{\mathbf{a}}-B_{t_{k}^{N}%
}^{\mathbf{a}})^{2}=(B_{t}^{\mathbf{a}})^{2}-2\int_{0}^{t}B_{s}^{\mathbf{a}%
}dB_{s}^{\mathbf{a}}. \label{quadra-def}%
\end{equation}
By the above construction, $\left \langle B^{\mathbf{a}}\right \rangle _{t}$,
$t\geq0$, is an increasing process with $\left \langle B^{\mathbf{a}%
}\right \rangle _{0}=0$. We call it the \textbf{quadratic variation process} of
the $G$--Brownian motion $B^{\mathbf{a}}$. Clearly $\left \langle
B^{\mathbf{a}}\right \rangle $ is an increasing process. It is also clear that,
for each $0\leq s\leq t$ and for each smooth real function $\psi$ such that
$\psi(\left \langle B^{\mathbf{a}}\right \rangle _{t-s})\in L_{G}^{1}%
(\mathcal{H}_{t-s})$ we have $\mathbb{E}[\psi(\left \langle B^{\mathbf{a}%
}\right \rangle _{t-s})]=\mathbb{E}[\psi(\left \langle B^{\mathbf{a}%
}\right \rangle _{t}-\left \langle B^{\mathbf{a}}\right \rangle _{s})]$. We also
have%
\[
\left \langle B^{\mathbf{a}}\right \rangle _{t}=\left \langle B^{-\mathbf{a}%
}\right \rangle _{t}=\left \langle -B^{\mathbf{a}}\right \rangle _{t}.
\]
It is important to keep in mind that $\left \langle B^{\mathbf{a}}\right \rangle
_{t}$ is not a deterministic process unless the case $\sigma_{\mathbf{aa}^{T}%
}=-\sigma_{-\mathbf{aa}^{T}}$ and thus $B^{\mathbf{a}}$ becomes a classical
Brownian motion. In fact we have

\begin{lemma}
\label{Lem-Q1}For each $0\leq s\leq t<\infty$%
\begin{align}
\mathbb{E}[\left \langle B^{\mathbf{a}}\right \rangle _{t}-\left \langle
B^{\mathbf{a}}\right \rangle _{s}|\mathcal{H}_{s}]  &  =\sigma_{\mathbf{aa}%
^{T}}(t-s),\  \  \label{quadra}\\
\mathbb{E}[-(\left \langle B^{\mathbf{a}}\right \rangle _{t}-\left \langle
B^{\mathbf{a}}\right \rangle _{s})|\mathcal{H}_{s}]  &  =\sigma_{-\mathbf{aa}%
^{T}}(t-s). \label{quadra1}%
\end{align}

\end{lemma}

\begin{proof}
By the definition of $\left \langle B^{\mathbf{a}}\right \rangle $ and
Proposition \ref{Prop-Integ}-(iii), then Example \ref{Exam-B2},
\begin{align*}
\mathbb{E}[\left \langle B^{\mathbf{a}}\right \rangle _{t}-\left \langle
B^{\mathbf{a}}\right \rangle _{s}|\mathcal{H}_{s}]  &  =\mathbb{E}%
[(B_{t}^{\mathbf{a}})^{2}-(B_{s}^{\mathbf{a}})^{2}-2\int_{s}^{t}%
B_{u}^{\mathbf{a}}dB_{u}^{\mathbf{a}}|\mathcal{H}_{s}]\\
&  =\mathbb{E}[(B_{t}^{\mathbf{a}})^{2}-(B_{s}^{\mathbf{a}})^{2}%
|\mathcal{H}_{s}]=\sigma_{\mathbf{aa}^{T}}(t-s).
\end{align*}
We then have (\ref{quadra}). (\ref{quadra1}) can be proved analogously by
using the equality $\mathbb{E}[-((B_{t}^{\mathbf{a}})^{2}-(B_{s}^{\mathbf{a}%
})^{2})|\mathcal{H}_{s}]=\sigma_{-\mathbf{aa}^{T}}(t-s)$.
\end{proof}

An interesting new phenomenon of our $G$-Brownian motion is that its quadratic
process $\left \langle B\right \rangle $ also has independent increments. In
fact, we have

\begin{lemma}
An increment of $\left \langle B^{\mathbf{a}}\right \rangle $ is the quadratic
variation of the corresponding increment of $B^{\mathbf{a}}$, i.e., for each
fixed $s\geq0$,%
\[
\left \langle B^{\mathbf{a}}\right \rangle _{t+s}-\left \langle B^{\mathbf{a}%
}\right \rangle _{s}=\left \langle (B^{s})^{\mathbf{a}}\right \rangle _{t}%
\]
where $B_{t}^{s}=B_{t+s}-B_{s}$, $t\geq0$ and $(B^{s})_{t}^{\mathbf{a}%
}=(\mathbf{a},B_{s}^{t})$.
\end{lemma}

\begin{proof}%
\begin{align*}
\left \langle B^{\mathbf{a}}\right \rangle _{t+s}-\left \langle B^{\mathbf{a}%
}\right \rangle _{s}  &  =(B_{t+s}^{\mathbf{a}})^{2}-2\int_{0}^{t+s}%
B_{u}^{\mathbf{a}}dB_{u}^{\mathbf{a}}-\left(  (B_{s}^{\mathbf{a}})^{2}%
-2\int_{0}^{s}B_{u}^{\mathbf{a}}dB_{u}^{\mathbf{a}}\right) \\
&  =(B_{t+s}^{\mathbf{a}}-B_{s}^{\mathbf{a}})^{2}-2\int_{s}^{t+s}%
(B_{u}^{\mathbf{a}}-B_{s}^{\mathbf{a}})dB_{u}^{\mathbf{a}}\\
&  =(B_{t+s}^{\mathbf{a}}-B_{s}^{\mathbf{a}})^{2}-2\int_{s}^{t}(B_{s+u}%
^{\mathbf{a}}-B_{s}^{\mathbf{a}})d(B_{u}^{\mathbf{a}}-B_{s}^{\mathbf{a}})\\
&  =\left \langle (B^{s})^{\mathbf{a}}\right \rangle _{t}.
\end{align*}

\end{proof}

\begin{lemma}
We have
\begin{equation}
\mathbb{E}[\left \langle B^{\mathbf{a}}\right \rangle _{t}^{2}]=\mathbb{E}%
[(\left \langle B^{\mathbf{a}}\right \rangle _{t+s}-\left \langle B^{\mathbf{a}%
}\right \rangle _{s})^{2}|\mathcal{H}_{s}]=\sigma_{\mathbf{aa}^{T}}^{2}%
t^{2},\  \ s,t\geq0. \label{Qua2}%
\end{equation}

\end{lemma}

\begin{proof}
We set $\phi(t):=\mathbb{E}[\left \langle B^{\mathbf{a}}\right \rangle _{t}%
^{2}]$.
\begin{align*}
\phi(t)  &  =\mathbb{E}[\{(B_{t}^{\mathbf{a}})^{2}-2\int_{0}^{t}%
B_{u}^{\mathbf{a}}dB_{u}^{\mathbf{a}}\}^{2}]\\
&  \leq2\mathbb{E}[(B_{t}^{\mathbf{a}})^{4}]+8\mathbb{E}[(\int_{0}^{t}%
B_{u}^{\mathbf{a}}dB_{u}^{\mathbf{a}})^{2}]\\
&  \leq6\sigma_{\mathbf{aa}^{T}}^{2}t^{2}+8\sigma_{\mathbf{aa}^{T}}\int
_{0}^{t}\mathbb{E[(}B_{u}^{\mathbf{a}})^{2}]du\\
&  =10\sigma_{\mathbf{aa}^{T}}^{2}t^{2}.
\end{align*}
This also implies $\mathbb{E}[(\left \langle B^{\mathbf{a}}\right \rangle
_{t}-\left \langle B^{\mathbf{a}}\right \rangle _{s})^{2}]=\phi(t-s)\leq
10\sigma_{\mathbf{aa}^{T}}^{2}(t-s)^{2}$. For each $s\in \lbrack0,t)$,
\begin{align*}
\phi(t)  &  =\mathbb{E}[(\left \langle B^{\mathbf{a}}\right \rangle
_{s}+\left \langle B^{\mathbf{a}}\right \rangle _{t}-\left \langle B^{\mathbf{a}%
}\right \rangle _{s})^{2}]\\
&  \leq \mathbb{E}[(\left \langle B^{\mathbf{a}}\right \rangle _{s}%
)^{2}]+\mathbb{E}[(\left \langle B^{\mathbf{a}}\right \rangle _{t}-\left \langle
B^{\mathbf{a}}\right \rangle _{s})^{2}]+2\mathbb{E}[(\left \langle
B^{\mathbf{a}}\right \rangle _{t}-\left \langle B^{\mathbf{a}}\right \rangle
_{s})\left \langle B^{\mathbf{a}}\right \rangle _{s}]\\
&  =\phi(s)+\phi(t-s)+2\mathbb{E}[\mathbb{E}[(\left \langle B^{\mathbf{a}%
}\right \rangle _{t}-\left \langle B^{\mathbf{a}}\right \rangle _{s}%
)|\mathcal{H}_{s}]\left \langle B^{\mathbf{a}}\right \rangle _{s}]\\
&  =\phi(s)+\phi(t-s)+2\sigma_{\mathbf{aa}^{T}}^{2}s(t-s).
\end{align*}
We set $\delta_{N}=t/N$, $t_{k}^{N}=kt/N=k\delta_{N}$ for a positive integer
$N$. By the above inequalities%
\begin{align*}
\phi(t_{N}^{N})  &  \leq \phi(t_{N-1}^{N})+\phi(\delta_{N})+2\sigma
_{\mathbf{aa}^{T}}^{2}t_{N-1}^{N}\delta_{N}\\
&  \leq \phi(t_{N-2}^{N})+2\phi(\delta_{N})+2\sigma_{\mathbf{aa}^{T}}%
^{2}(t_{N-1}^{N}+t_{N-2}^{N})\delta_{N}\\
&  \vdots \\
&
\end{align*}
We then have
\[
\phi(t)\leq N\phi(\delta_{N})+2\sigma_{\mathbf{aa}^{T}}^{2}\sum_{k=0}%
^{N-1}t_{k}^{N}\delta_{N}\leq14t^{2}\sigma_{\mathbf{aa}^{T}}^{2}%
/N+2\sigma_{\mathbf{aa}^{T}}^{2}\sum_{k=0}^{N-1}t_{k}^{N}\delta_{N}.
\]
Let $N\rightarrow \infty$ we have $\phi(t)\leq2\sigma_{\mathbf{aa}^{T}}^{2}%
\int_{0}^{t}sds=\sigma_{\mathbf{aa}^{T}}^{2}t^{2}$. Thus $\mathbb{E}%
[\left \langle B^{\mathbf{a}}\right \rangle _{t}^{2}]\leq \sigma_{\mathbf{aa}%
^{T}}^{2}t^{2}$. This with $\mathbb{E}[\left \langle B^{\mathbf{a}%
}\right \rangle _{t}^{2}]\geq E^{0}[\left \langle B^{\mathbf{a}}\right \rangle
_{t}^{2}]=\sigma_{\mathbf{aa}^{T}}^{2}t^{2}$ implies (\ref{Qua2}). In the last
step, the classical normal distribution $P_{1}^{0}$, or $N(0,\gamma_{0}%
\gamma_{0}^{T})$, $\gamma_{0}\in \Gamma$, is chosen such that%
\[
tr[\gamma_{0}\gamma_{0}^{T}\mathbf{aa}^{T}]=\sigma_{\mathbf{aa}^{T}}^{2}%
=\sup_{\gamma \in \Gamma}tr[\gamma \gamma^{T}\mathbf{aa}^{T}].
\]

\end{proof}

Similarly we have
\begin{align}
\mathbb{E}[(\left \langle B^{\mathbf{a}}\right \rangle _{t}-\left \langle
B^{\mathbf{a}}\right \rangle _{s})^{3}|\mathcal{H}_{s}]  &  =\sigma
_{\mathbf{aa}^{T}}^{3}(t-s)^{3},\label{Quad3}\\
\mathbb{E}[(\left \langle B^{\mathbf{a}}\right \rangle _{t}-\left \langle
B^{\mathbf{a}}\right \rangle _{s})^{4}|\mathcal{H}_{s}]  &  =\sigma
_{\mathbf{aa}^{T}}^{4}(t-s)^{4}.\nonumber
\end{align}

\begin{proposition}
\label{Prop-temp}Let $0\leq s\leq t$, $\xi \in L_{G}^{1}(\mathcal{H}_{s})$,
$X\in L_{G}^{1}(\mathcal{H})$. Then%
\begin{align*}
\mathbb{E}[X+\xi((B_{t}^{\mathbf{a}})^{2}-(B_{s}^{\mathbf{a}})^{2})]  &
=\mathbb{E}[X+\xi(B_{t}^{\mathbf{a}}-B_{s}^{\mathbf{a}})^{2}]\\
&  =\mathbb{E}[X+\xi(\left \langle B^{\mathbf{a}}\right \rangle _{t}%
-\left \langle B^{\mathbf{a}}\right \rangle _{s})].
\end{align*}

\end{proposition}

\begin{proof}
By (\ref{quadra-def}) and applying Proposition \ref{E-x+y}, we have%
\begin{align*}
\mathbb{E}[X+\xi((B_{t}^{\mathbf{a}})^{2}-(B_{s}^{\mathbf{a}})^{2})]  &
=\mathbb{E}[X+\xi(\left \langle B^{\mathbf{a}}\right \rangle _{t}-\left \langle
B^{\mathbf{a}}\right \rangle _{s}+2\int_{s}^{t}B_{u}^{\mathbf{a}}%
dB_{u}^{\mathbf{a}})]\\
&  =\mathbb{E}[X+\xi(\left \langle B^{\mathbf{a}}\right \rangle _{t}%
-\left \langle B^{\mathbf{a}}\right \rangle _{s})].
\end{align*}
We also have
\begin{align*}
\mathbb{E}[X+\xi((B_{t}^{\mathbf{a}})^{2}-(B_{s}^{\mathbf{a}})^{2})]  &
=\mathbb{E}[X+\xi \{(B_{t}^{\mathbf{a}}-B_{s}^{\mathbf{a}})^{2}+2(B_{t}%
^{\mathbf{a}}-B_{s}^{\mathbf{a}})B_{s}^{\mathbf{a}}\}]\\
&  =\mathbb{E}[X+\xi(B_{t}^{\mathbf{a}}-B_{s}^{\mathbf{a}})^{2}].
\end{align*}

\end{proof}

\begin{example}
\label{Exa-AB}We assume that in a financial market a stock price
$(S_{t})_{t\geq0}$ is observed. Let $B_{t}=\log(S_{t})$, $t\geq0$, be a
$1$-dimensional $G$-Brownian motion $(d=1)$ with $\Gamma=[\sigma_{\ast}%
,\sigma^{\ast}]$, with fixed $\sigma_{\ast}\in \lbrack0,\frac{1}{2})$ and
$\sigma^{\ast}\in \lbrack1,\infty)$. Two traders $a$ and $b$ in a same bank are
using their own statistics to price a contingent claim $X=\left \langle
B\right \rangle _{T}$ with maturity $T$. Suppose, for example, under the
probability measure $\mathbb{P}_{a}$ of $a$, $B$ is a (classical) Brownian
motion whereas under $\mathbb{P}_{b}$ of $b$, $\frac{1}{2}B$ is a Brownian
motion, here $\mathbb{P}_{a}$ (resp. $\mathbb{P}_{b}$) is a classical
probability measure with its linear expectation $\mathbb{E}^{a}$ (resp.
$\mathbb{E}^{b}$ ) generated by the heat equation $\partial_{t}u=\frac{1}%
{2}\partial_{xx}^{2}u$ (resp. $\partial_{t}u=\frac{1}{4}\partial_{xx}^{2}u$).
Since $\mathbb{E}^{a}$ and $\mathbb{E}^{b}$ are both dominated by $\mathbb{E}$
in the sense of (3), they can be both well--defined as a linear bounded
functional in $L_{G}^{1}(\mathcal{H})$. This framework cannot be provided by
just using a classical probability space because it is known that
$\left \langle B\right \rangle _{T}=T$, $\mathbb{P}^{a}$--a.s., and
$\left \langle B\right \rangle _{T}=\frac{T}{4}$, $\mathbb{P}^{b}$--a.s. Thus
there is no probability measure on $\Omega$ with respect to which $P_{a}$ and
$P_{b}$ are both absolutely continuous. Practically this sublinear expectation
$\mathbb{E}$ provides a realistic tool of dynamic risk measure for a risk
supervisor of the traders $a$ and $b$: given a risk position $X\in L_{G}%
^{1}(\mathcal{H}_{T})$ we always have $\mathbb{E}[-X|\mathcal{H}_{t}%
]\geq \mathbb{E}^{a}[-X|\mathcal{H}_{t}]\vee \mathbb{E}^{b}[-X|\mathcal{H}_{t}]$
for the loss $-X$ of this position. The meaning is that the supervisor uses a
more sensitive risk measure. Clearly no linear expectation can play this role.
The subset $\Gamma$ represents the uncertainty of the volatility model of a
risk regulator. The lager the subset $\Gamma$, the bigger the uncertainty,
thus the stronger the corresponding $\mathbb{E}$. \newline It is worth to
consider to create a hierarchic and dynamic risk control system for a bank, or
a banking system, in which the Chief Risk Officer (CRO) uses $\mathbb{E}%
=\mathbb{E}^{G}$ for his risk measure and the Risk Officer the $i$th division
of the bank uses $\mathbb{E}^{i}=\mathbb{E}^{G_{i}}$ for his one, where
\[
G(A)=\frac{1}{2}\sup_{\gamma \in \Gamma}\text{tr}[\gamma \gamma^{T}%
A],\  \ G_{i}(A)=\frac{1}{2}\sup_{\gamma \in \Gamma_{i}}\text{tr}[\gamma
\gamma^{T}A],\  \Gamma_{i}\subset \Gamma \text{,\ }\ i=1,\cdots,I.
\]
Thus $\mathbb{E}^{i}$ is dominated by $\mathbb{E}$ for each $i$. For a large
banking system we can even consider to create $\mathbb{E}^{ij}=\mathbb{E}%
^{G_{ij}}$ for its $(i,j)$th sub-division. The reasoning is: in general, a
risk regulator's statistics and knowledge of a specific risk position $X$ are
less than a trader who is concretely involved in the business of the product
$X$.
\end{example}

To define the integration of a process $\eta \in M_{G}^{1}(0,T)$ with respect
to $d\left \langle B^{\mathbf{a}}\right \rangle $, we first define a mapping:%
\[
Q_{0,T}(\eta)=\int_{0}^{T}\eta(s)d\left \langle B^{\mathbf{a}}\right \rangle
_{s}:=\sum_{k=0}^{N-1}\xi_{k}(\left \langle B^{\mathbf{a}}\right \rangle
_{t_{k+1}}-\left \langle B^{\mathbf{a}}\right \rangle _{t_{k}}):M_{G}%
^{1,0}(0,T)\mapsto L^{1}(\mathcal{H}_{T}).
\]

\begin{lemma}
\label{Lem-Q2}For each $\eta \in M_{G}^{1,0}(0,T)$,
\begin{equation}
\mathbb{E}[|Q_{0,T}(\eta)|]\leq \sigma_{\mathbf{aa}^{T}}\int_{0}^{T}%
\mathbb{E}[|\eta_{s}|]ds,\  \label{dA}%
\end{equation}
Thus $Q_{0,T}:M_{G}^{1,0}(0,T)\mapsto L^{1}(\mathcal{H}_{T})$ is a continuous
linear mapping. Consequently, $Q_{0,T}$ can be uniquely extended to $M_{G}%
^{1}(0,T)$. We still denote this mapping \ by%
\[
\int_{0}^{T}\eta(s)d\left \langle B^{\mathbf{a}}\right \rangle _{s}=Q_{0,T}%
(\eta),\  \  \eta \in M_{G}^{1}(0,T)\text{.}%
\]
We still have
\begin{equation}
\mathbb{E}[|\int_{0}^{T}\eta(s)d\left \langle B^{\mathbf{a}}\right \rangle
_{s}|]\leq \sigma_{\mathbf{aa}^{T}}\int_{0}^{T}\mathbb{E}[|\eta_{s}%
|]ds,\  \  \forall \eta \in M_{G}^{1}(0,T)\text{.} \label{qua-ine}%
\end{equation}

\end{lemma}

\begin{proof}
By applying Lemma \ref{Lem-Q1}, (\ref{dA}) can be checked as follows:%
\begin{align*}
\mathbb{E}[|\sum_{k=0}^{N-1}\xi_{k}(\left \langle B^{\mathbf{a}}\right \rangle
_{t_{k+1}}-\left \langle B^{\mathbf{a}}\right \rangle _{t_{k}})|]  &  \leq
\sum_{k=0}^{N-1}\mathbb{E[}|\xi_{k}|\cdot \mathbb{E}[\left \langle
B^{\mathbf{a}}\right \rangle _{t_{k+1}}-\left \langle B^{\mathbf{a}%
}\right \rangle _{t_{k}}|\mathcal{H}_{t_{k}}]]\\
&  =\sum_{k=0}^{N-1}\mathbb{E[}|\xi_{k}|]\sigma_{\mathbf{aa}^{T}}%
(t_{k+1}-t_{k})\\
&  =\sigma_{\mathbf{aa}^{T}}\int_{0}^{T}\mathbb{E}[|\eta_{s}|]ds.
\end{align*}

\end{proof}

We have the following isometry

\begin{proposition}
Let $\eta \in M_{G}^{2}(0,T)$%
\begin{equation}
\mathbb{E}[(\int_{0}^{T}\eta(s)dB_{s}^{\mathbf{a}})^{2}]=\mathbb{E}[\int
_{0}^{T}\eta^{2}(s)d\left \langle B^{\mathbf{a}}\right \rangle _{s}]
\label{isometry}%
\end{equation}

\end{proposition}

\begin{proof}
We first consider $\eta \in M_{G}^{2,0}(0,T)$ with the form
\[
\eta_{t}(\omega)=\sum_{k=0}^{N-1}\xi_{k}(\omega)\mathbf{I}_{[t_{k},t_{k+1}%
)}(t)
\]
and thus $\int_{0}^{T}\eta(s)dB_{s}^{\mathbf{a}}:=\sum_{k=0}^{N-1}\xi
_{k}(B_{t_{k+1}}^{\mathbf{a}}-B_{t_{k}}^{\mathbf{a}})$\textbf{.} By
Proposition \ref{E-x+y} we have
\[
\mathbb{E}[X+2\xi_{k}(B_{t_{k+1}}^{\mathbf{a}}-B_{t_{k}}^{\mathbf{a}})\xi
_{l}(B_{t_{l+1}}^{\mathbf{a}}-B_{t_{l}}^{\mathbf{a}})]=\mathbb{E}[X]\text{,
for }X\in L_{G}^{1}(\mathcal{H)}\text{, }l\not =k.
\]
Thus%
\[
\mathbb{E}[(\int_{0}^{T}\eta(s)dB_{s}^{\mathbf{a}})^{2}]=\mathbb{E[}\left(
\sum_{k=0}^{N-1}\xi_{k}(B_{t_{k+1}}^{\mathbf{a}}-B_{t_{k}}^{\mathbf{a}%
})\right)  ^{2}]=\mathbb{E[}\sum_{k=0}^{N-1}\xi_{k}^{2}(B_{t_{k+1}%
}^{\mathbf{a}}-B_{t_{k}}^{\mathbf{a}})^{2}].
\]
This with Proposition \ref{Prop-temp}, it follows that
\[
\mathbb{E}[(\int_{0}^{T}\eta(s)dB_{s}^{\mathbf{a}})^{2}]=\mathbb{E[}\sum
_{k=0}^{N-1}\xi_{k}^{2}(\left \langle B^{\mathbf{a}}\right \rangle _{t_{k+1}%
}-\left \langle B^{\mathbf{a}}\right \rangle _{t_{k}})]=\mathbb{E[}\int_{0}%
^{T}\eta^{2}(s)d\left \langle B^{\mathbf{a}}\right \rangle _{s}].
\]
Thus (\ref{isometry}) holds for $\eta \in M_{G}^{2,0}(0,T)$. We thus can
continuously extend this equality to the case $\eta \in M_{G}^{2}(0,T)$ and
obtain (\ref{isometry}).
\end{proof}

\subsection{Mutual variation processes for $G$--Brownian motion}

Let $\mathbf{a}=(a_{1},\cdots,a_{d})^{T}$ and $\mathbf{\bar{a}}=(\bar{a}%
_{1},\cdots,\bar{a}_{d})^{T}$ be two given vectors in $\mathbb{R}^{d}$. We
then have their quadratic variation process $\left \langle B^{\mathbf{a}%
}\right \rangle $ and $\left \langle B^{\mathbf{\bar{a}}}\right \rangle $. We
then can define their mutual variation process by%
\begin{align*}
\left \langle B^{\mathbf{a}},B^{\mathbf{\bar{a}}}\right \rangle _{t}  &
:=\frac{1}{4}[\left \langle B^{\mathbf{a}}+B^{\mathbf{\bar{a}}}\right \rangle
_{t}-\left \langle B^{\mathbf{a}}-B^{\mathbf{\bar{a}}}\right \rangle _{t}]\\
&  =\frac{1}{4}[\left \langle B^{\mathbf{a}+\mathbf{\bar{a}}}\right \rangle
_{t}-\left \langle B^{\mathbf{a}-\mathbf{\bar{a}}}\right \rangle _{t}].
\end{align*}
Since $\left \langle B^{\mathbf{a}-\mathbf{\bar{a}}}\right \rangle =\left \langle
B^{\mathbf{\bar{a}}-\mathbf{a}}\right \rangle =\left \langle -B^{\mathbf{a}%
-\mathbf{\bar{a}}}\right \rangle $, we see that $\left \langle B^{\mathbf{a}%
},B^{\mathbf{\bar{a}}}\right \rangle _{t}=\left \langle B^{\mathbf{\bar{a}}%
},B^{\mathbf{a}}\right \rangle _{t}$. In particular we have $\left \langle
B^{\mathbf{a}},B^{\mathbf{a}}\right \rangle =\left \langle B^{\mathbf{a}%
}\right \rangle $. Let $\pi_{t}^{N}$, $N=1,2,\cdots$, be a sequence of
partitions of $[0,t]$. We observe that%
\[
\sum_{k=0}^{N-1}(B_{t_{k+1}^{N}}^{\mathbf{a}}-B_{t_{k}^{N}}^{\mathbf{a}%
})(B_{t_{k+1}^{N}}^{\mathbf{\bar{a}}}-B_{t_{k}^{N}}^{\mathbf{\bar{a}}}%
)=\frac{1}{4}\sum_{k=0}^{N-1}[(B_{t_{k+1}}^{\mathbf{a}+\mathbf{\bar{a}}%
}-B_{t_{k}}^{\mathbf{a}+\mathbf{\bar{a}}})^{2}-(B_{t_{k+1}}^{\mathbf{a}%
-\mathbf{\bar{a}}}-B_{t_{k}}^{\mathbf{a}-\mathbf{\bar{a}}})^{2}].
\]
Thus as $\mu(\pi_{t}^{N})\rightarrow0$, we have%
\[
\lim_{N\rightarrow0}\sum_{k=0}^{N-1}(B_{t_{k+1}^{N}}^{\mathbf{a}}-B_{t_{k}%
^{N}}^{\mathbf{a}})(B_{t_{k+1}^{N}}^{\mathbf{\bar{a}}}-B_{t_{k}^{N}%
}^{\mathbf{\bar{a}}})=\left \langle B^{\mathbf{a}},B^{\mathbf{\bar{a}}%
}\right \rangle _{t}.
\]
We also have
\begin{align*}
\left \langle B^{\mathbf{a}},B^{\mathbf{\bar{a}}}\right \rangle _{t}  &
=\frac{1}{4}[\left \langle B^{\mathbf{a}+\mathbf{\bar{a}}}\right \rangle
_{t}-\left \langle B^{\mathbf{a}-\mathbf{\bar{a}}}\right \rangle _{t}]\\
&  =\frac{1}{4}[(B_{t}^{\mathbf{a}+\mathbf{\bar{a}}})^{2}-2\int_{0}^{t}%
B_{s}^{\mathbf{a}+\mathbf{\bar{a}}}dB_{s}^{\mathbf{a}+\mathbf{\bar{a}}}%
-(B_{t}^{\mathbf{a}-\mathbf{\bar{a}}})^{2}+2\int_{0}^{t}B_{s}^{\mathbf{a}%
-\mathbf{\bar{a}}}dB_{s}^{\mathbf{a}-\mathbf{\bar{a}}}]\\
&  =B_{t}^{\mathbf{a}}B_{t}^{\mathbf{\bar{a}}}-\int_{0}^{t}B_{s}^{\mathbf{a}%
}dB_{s}^{\mathbf{\bar{a}}}-\int_{0}^{t}B_{s}^{\mathbf{\bar{a}}}dB_{s}%
^{\mathbf{a}}.
\end{align*}
Now for each $\eta \in M_{G}^{1}(0,T)$ we can consistently define
\[
\int_{0}^{T}\eta_{s}d\left \langle B^{\mathbf{a}},B^{\mathbf{\bar{a}}%
}\right \rangle _{s}=\frac{1}{4}\int_{0}^{T}\eta_{s}d\left \langle
B^{\mathbf{a}+\mathbf{\bar{a}}}\right \rangle _{s}-\frac{1}{4}\int_{0}^{T}%
\eta_{s}\left \langle B^{\mathbf{a}-\mathbf{\bar{a}}}\right \rangle _{s}.
\]

\begin{lemma}
\label{Lem-mutual}Let $\eta^{N}\in M_{G}^{1,0}(0,T)$, $N=1,2,\cdots,$ be of
form
\[
\eta_{t}^{N}(\omega)=\sum_{k=0}^{N-1}\xi_{k}^{N}(\omega)\mathbf{I}_{[t_{k}%
^{N},t_{k+1}^{N})}(t)
\]
with $\mu(\pi_{T}^{N})\rightarrow0$ and $\eta^{N}\rightarrow \eta$ in
$M_{G}^{1}(0,T)$ as $N\rightarrow \infty$. Then we have the following
convergence in $L_{G}^{1}(\mathcal{H}_{T})$:%
\begin{align*}
\int_{0}^{T}\eta^{N}(s)d\left \langle B^{\mathbf{a}},B^{\mathbf{\bar{a}}%
}\right \rangle _{s}  &  :=\sum_{k=0}^{N-1}\xi_{k}^{N}(B_{t_{k+1}^{N}%
}^{\mathbf{a}}-B_{t_{k}^{N}}^{\mathbf{a}})(B_{t_{k+1}^{N}}^{\mathbf{\bar{a}}%
}-B_{t_{k}^{N}}^{\mathbf{\bar{a}}})\\
&  \rightarrow \int_{0}^{T}\eta(s)d\left \langle B^{\mathbf{a}},B^{\mathbf{\bar
{a}}}\right \rangle _{s}.
\end{align*}

\end{lemma}

\subsection{It\^{o}'s formula for $G$--Brownian motion}

We have the corresponding It\^{o}'s formula of $\Phi(X_{t})$ for a
\textquotedblleft$G$-It\^{o} process\textquotedblright \ $X$. For
simplification, we only treat the case where the function $\Phi$ is
sufficiently regular. For notational simplification, we denote $B^{i}%
=B^{\mathbf{e}_{i}}$, the $i$--th coordinate of the $G$--Brownian motion $B$,
under a given orthonormal base $(\mathbf{e}_{1},\cdots,\mathbf{e}_{d})$ of
$\mathbb{R}^{d}$.

\begin{lemma}
\label{Lem-26}Let $\Phi \in C^{2}(\mathbb{R}^{n})$ be bounded with bounded
derivatives and $\{ \partial_{x^{\mu}x^{\nu}}^{2}\Phi \}_{\mu,\nu=1}^{n}$ are
uniformly Lipschitz. Let $s\in \lbrack0,T]$ be fixed and let $X=(X^{1}%
,\cdots,X^{n})^{T}$ be an $n$--dimensional process on $[s,T]$ of the form
\[
X_{t}^{\nu}=X_{s}^{\nu}+\alpha^{\nu}(t-s)+\eta^{\nu ij}(\left \langle
B^{i},B^{j}\right \rangle _{t}-\left \langle B^{i},B^{j}\right \rangle
_{s})+\beta^{\nu j}(B_{t}^{j}-B_{s}^{j}),
\]
where, for $\nu=1,\cdots,n$, $i,j=1,\cdots,d$, $\alpha^{\nu}$, $\eta^{\nu ij}$
and $\beta^{\nu ij}$, are bounded elements of $L_{G}^{2}(\mathcal{H}_{s})$ and
$X_{s}=(X_{s}^{1},\cdots,X_{s}^{n})^{T}$ is a given $\mathbb{R}^{n}$--vector
in $L_{G}^{2}(\mathcal{H}_{s})$. Then we have
\begin{align}
\Phi(X_{t})-\Phi(X_{s})  &  =\int_{s}^{t}\partial_{x^{\nu}}\Phi(X_{u}%
)\beta^{\nu j}dB_{u}^{j}+\int_{s}^{t}\partial_{x_{\nu}}\Phi(X_{u})\alpha^{\nu
}du\label{B-Ito}\\
&  +\int_{s}^{t}[\partial_{x^{\nu}}\Phi(X_{u})\eta^{\nu ij}+\frac{1}%
{2}\partial_{x^{\mu}x^{\nu}}^{2}\Phi(X_{u})\beta^{\nu i}\beta^{\nu
j}]d\left \langle B^{i},B^{j}\right \rangle _{u}.\nonumber
\end{align}
Here we use the Einstein convention, i.e., the above repeated indices $\mu
,\nu$, $i$ and $j$ (but not $k$) imply the summation.
\end{lemma}

\begin{proof}
For each positive integer $N$ we set $\delta=(t-s)/N$ and take the partition
\[
\pi_{\lbrack s,t]}^{N}=\{t_{0}^{N},t_{1}^{N},\cdots,t_{N}^{N}\}=\{s,s+\delta
,\cdots,s+N\delta=t\}.
\]
We have
\begin{align}
\Phi(X_{t})-\Phi(X_{s})  &  =\sum_{k=0}^{N-1}[\Phi(X_{t_{k+1}^{N}}%
)-\Phi(X_{t_{k}^{N}})]\nonumber \\
&  =\sum_{k=0}^{N-1}[\partial_{x^{\mu}}\Phi(X_{t_{k}^{N}})(X_{t_{k+1}^{N}%
}^{\mu}-X_{t_{k}^{N}}^{\mu})\nonumber \\
&  +\frac{1}{2}[\partial_{x^{\mu}x^{\nu}}^{2}\Phi(X_{t_{k}^{N}})(X_{t_{k+1}%
^{N}}^{\mu}-X_{t_{k}^{N}}^{\mu})(X_{t_{k+1}^{N}}^{\nu}-X_{t_{k}^{N}}^{\nu
})+\eta_{k}^{N}]] \label{Ito}%
\end{align}
where
\[
\eta_{k}^{N}=[\partial_{x^{\mu}x^{\nu}}^{2}\Phi(X_{t_{k}^{N}}+\theta
_{k}(X_{t_{k+1}^{N}}-X_{t_{k}^{N}}))-\partial_{x^{\mu}x^{\nu}}^{2}%
\Phi(X_{t_{k}^{N}})](X_{t_{k+1}^{N}}^{\mu}-X_{t_{k}^{N}}^{\mu})(X_{t_{k+1}%
^{N}}^{\nu}-X_{t_{k}^{N}}^{\nu})
\]
with $\theta_{k}\in \lbrack0,1]$. We have%
\begin{align*}
\mathbb{E}[|\eta_{k}^{N}|]  &  =\mathbb{E}[|[\partial_{x^{\mu}x^{\nu}}^{2}%
\Phi(X_{t_{k}^{N}}+\theta_{k}(X_{t_{k+1}^{N}}-X_{t_{k}^{N}}))-\partial
_{x^{\mu}x^{\nu}}^{2}\Phi(X_{t_{k}^{N}})]\\
&  \times(X_{t_{k+1}^{N}}^{\mu}-X_{t_{k}^{N}}^{\mu})(X_{t_{k+1}^{N}}^{\nu
}-X_{t_{k}^{N}}^{\nu})|]\\
&  \leq c\mathbb{E[}|X_{t_{k+1}^{N}}-X_{t_{k}^{N}}|^{3}]\leq C[\delta
^{3}+\delta^{3/2}]
\end{align*}
where $c$ is the Lipschitz constant of $\{ \partial_{x^{\mu}x^{\nu}}^{2}%
\Phi \}_{\mu,\nu=1}^{d}$. In the last step we use Example \ref{Exam-1} and
(\ref{Quad3}). Thus $\sum_{k}\mathbb{E}[|\eta_{k}^{N}|]\rightarrow0$. The rest
terms in the summation of the right side of (\ref{Ito}) are $\xi_{t}^{N}%
+\zeta_{t}^{N}$ with%
\begin{align*}
\xi_{t}^{N}  &  =\sum_{k=0}^{N-1}\{ \partial_{x^{\mu}}\Phi(X_{t_{k}^{N}%
})[\alpha^{\mu}(t_{k+1}^{N}-t_{k}^{N})+\eta^{\mu ij}(\left \langle B^{i}%
,B^{j}\right \rangle _{t_{k+1}^{N}}-\left \langle B^{i},B^{j}\right \rangle
_{t_{k}^{N}})\\
&  +\beta^{\mu j}(B_{t_{k+1}^{N}}^{j}-B_{t_{k}^{N}}^{j})]+\frac{1}{2}%
\partial_{x^{\mu}x^{\nu}}^{2}\Phi(X_{t_{k}^{N}})\beta^{\mu i}\beta^{\nu
j}(B_{t_{k+1}^{N}}^{i}-B_{t_{k}^{N}}^{i})(B_{t_{k+1}^{N}}^{j}-B_{t_{k}^{N}%
}^{j})\} \\
&
\end{align*}
and
\begin{align*}
\zeta_{t}^{N}  &  =\frac{1}{2}\sum_{k=0}^{N-1}\partial_{x^{\mu}x^{\nu}}%
^{2}\Phi(X_{t_{k}^{N}})[\alpha^{\mu}(t_{k+1}^{N}-t_{k}^{N})+\eta^{\mu
ij}(\left \langle B^{i},B^{j}\right \rangle _{t_{k+1}^{N}}-\left \langle
B^{i},B^{j}\right \rangle _{t_{k}^{N}})]\\
&  \times \lbrack \alpha^{\nu}(t_{k+1}^{N}-t_{k}^{N})+\eta^{\nu lm}(\left \langle
B^{l},B^{m}\right \rangle _{t_{k+1}^{N}}-\left \langle B^{l},B^{m}\right \rangle
_{t_{k}^{N}})]\\
&  +[\alpha^{\mu}(t_{k+1}^{N}-t_{k}^{N})+\eta^{\mu ij}(\left \langle
B^{i},B^{j}\right \rangle _{t_{k+1}^{N}}-\left \langle B^{i},B^{j}\right \rangle
_{t_{k}^{N}})]\beta^{\nu l}(B_{t_{k+1}^{N}}^{l}-B_{t_{k}^{N}}^{l})
\end{align*}
We observe that, for each $u\in \lbrack t_{k}^{N},t_{k+1}^{N})$
\begin{align*}
&  \mathbb{E}[|\partial_{x^{\mu}}\Phi(X_{u})-\sum_{k=0}^{N-1}\partial_{x^{\mu
}}\Phi(X_{t_{k}^{N}})\mathbf{I}_{[t_{k}^{N},t_{k+1}^{N})}(u)|^{2}]\\
&  =\mathbb{E}[|\partial_{x^{\mu}}\Phi(X_{u})-\partial_{x^{\mu}}\Phi
(X_{t_{k}^{N}})|^{2}]\\
&  \leq c^{2}\mathbb{E}[|X_{u}-X_{t_{k}^{N}}|^{2}]\leq C[\delta+\delta^{2}].
\end{align*}
Thus $\sum_{k=0}^{N-1}\partial_{x^{\mu}}\Phi(X_{t_{k}^{N}})\mathbf{I}%
_{[t_{k}^{N},t_{k+1}^{N})}(\cdot)$ tends to $\partial_{x^{\mu}}\Phi(X_{\cdot
})$ in $M_{G}^{2}(0,T)$. Similarly,
\[
\sum_{k=0}^{N-1}\partial_{x^{\mu}x^{\nu}}^{2}\Phi(X_{t_{k}^{N}})\mathbf{I}%
_{[t_{k}^{N},t_{k+1}^{N})}(\cdot)\rightarrow \partial_{x^{\mu}x^{\nu}}^{2}%
\Phi(X_{\cdot})\text{, in \ }M_{G}^{2}(0,T).
\]
Let $N\rightarrow \infty$, by Lemma \ref{Lem-mutual} as well as the definitions
of the integrations of $dt$, $dB_{t}$ and $d\left \langle B\right \rangle _{t}$,
the limit of $\xi_{t}^{N}$ in $L_{G}^{2}(\mathcal{H}_{t})$ is just the right
hand side of (\ref{B-Ito}). By the next Remark, we also have $\zeta_{t}%
^{N}\rightarrow0$ in $L_{G}^{2}(\mathcal{H}_{t})$. We then have proved
(\ref{B-Ito}).
\end{proof}

\begin{remark}
In the proof of $\zeta_{t}^{N}\rightarrow0$ in $L_{G}^{2}(\mathcal{H}_{t})$,
we use the following estimates: for $\psi^{N}\in M_{G}^{1,0}(0,T)$ such that
$\psi_{t}^{N}=\sum_{k=0}^{N-1}\xi_{t_{k}}^{N}\mathbf{I}_{[t_{k}^{N}%
,t_{k+1}^{N})}(t)$, and $\pi_{T}^{N}=\{0\leq t_{0},\cdots,t_{N}=T\}$ with
$\lim_{N\rightarrow \infty}\mu(\pi_{T}^{N})=0$ and $\sum_{k=0}^{N-1}%
\mathbb{E}[|\xi_{t_{k}}^{N}|](t_{k+1}^{N}-t_{k}^{N})\leq C$, for all
$N=1,2,\cdots$, we have $\mathbb{E}[|\sum_{k=0}^{N-1}\xi_{k}^{N}(t_{k+1}%
^{N}-t_{k}^{N})^{2}]\rightarrow0$ and, for any fixed $\mathbf{a,\bar{a}\in
}\mathbb{R}^{d}$,%
\begin{align*}
\mathbb{E}[|\sum_{k=0}^{N-1}\xi_{k}^{N}(\left \langle B^{\mathbf{a}%
}\right \rangle _{t_{k+1}^{N}}-\left \langle B^{\mathbf{a}}\right \rangle
_{t_{k}^{N}})^{2}|]  &  \leq \sum_{k=0}^{N-1}\mathbb{E[}|\xi_{k}^{N}%
|\cdot \mathbb{E}[(\left \langle B^{\mathbf{a}}\right \rangle _{t_{k+1}^{N}%
}-\left \langle B^{\mathbf{a}}\right \rangle _{t_{k}^{N}})^{2}|\mathcal{H}%
_{t_{k}^{N}}]]\\
&  =\sum_{k=0}^{N-1}\mathbb{E[}|\xi_{k}^{N}|]\sigma_{\mathbf{aa}^{T}}%
^{2}(t_{k+1}^{N}-t_{k}^{N})^{2}\rightarrow0,
\end{align*}%
\begin{align*}
&  \mathbb{E}[|\sum_{k=0}^{N-1}\xi_{k}^{N}(\left \langle B^{\mathbf{a}%
}\right \rangle _{t_{k+1}^{N}}-\left \langle B^{\mathbf{a}}\right \rangle
_{t_{k}^{N}})(t_{k+1}^{N}-t_{k}^{N})|]\\
&  \leq \sum_{k=0}^{N-1}\mathbb{E[}|\xi_{k}^{N}|(t_{k+1}^{N}-t_{k}^{N}%
)\cdot \mathbb{E}[(\left \langle B^{\mathbf{a}}\right \rangle _{t_{k+1}^{N}%
}-\left \langle B^{\mathbf{a}}\right \rangle _{t_{k}^{N}})|\mathcal{H}%
_{t_{k}^{N}}]]\\
&  =\sum_{k=0}^{N-1}\mathbb{E[}|\xi_{k}^{N}|]\sigma_{\mathbf{aa}^{T}}%
(t_{k+1}^{N}-t_{k}^{N})^{2}\rightarrow0,
\end{align*}
as well as
\begin{align*}
\mathbb{E}[|\sum_{k=0}^{N-1}\xi_{k}^{N}(t_{k+1}^{N}-t_{k}^{N})(B_{t_{k+1}^{N}%
}^{\mathbf{a}}-B_{t_{k}^{N}}^{\mathbf{a}})]|  &  \leq \sum_{k=0}^{N-1}%
\mathbb{E[}|\xi_{k}^{N}|](t_{k+1}^{N}-t_{k}^{N})\mathbb{E}[|B_{t_{k+1}^{N}%
}^{\mathbf{a}}-B_{t_{k}^{N}}^{\mathbf{a}}|]\\
&  =\sqrt{\frac{2\sigma_{\mathbf{aa}^{T}}}{\pi}}\sum_{k=0}^{N-1}%
\mathbb{E[}|\xi_{k}^{N}|](t_{k+1}^{N}-t_{k}^{N})^{3/2}\rightarrow0\
\end{align*}
and%
\begin{align*}
&  \mathbb{E}[|\sum_{k=0}^{N-1}\xi_{k}^{N}(\left \langle B^{\mathbf{a}%
}\right \rangle _{t_{k+1}^{N}}-\left \langle B^{\mathbf{a}}\right \rangle
_{t_{k}^{N}})(B_{t_{k+1}^{N}}^{\mathbf{\bar{a}}}-B_{t_{k}^{N}}^{\mathbf{\bar
{a}}})|]\\
&  \leq \sum_{k=0}^{N-1}\mathbb{E[}|\xi_{k}^{N}|]\mathbb{E[}(\left \langle
B^{\mathbf{a}}\right \rangle _{t_{k+1}^{N}}-\left \langle B^{\mathbf{a}%
}\right \rangle _{t_{k}^{N}})|B_{t_{k+1}^{N}}^{\mathbf{\bar{a}}}-B_{t_{k}^{N}%
}^{\mathbf{\bar{a}}}|]\\
&  \leq \sum_{k=0}^{N-1}\mathbb{E[}|\xi_{k}^{N}|]\mathbb{E[}(\left \langle
B^{\mathbf{a}}\right \rangle _{t_{k+1}^{N}}-\left \langle B^{\mathbf{a}%
}\right \rangle _{t_{k}^{N}})^{2}]^{1/2}\mathbb{E[}|B_{t_{k+1}^{N}%
}^{\mathbf{\bar{a}}}-B_{t_{k}^{N}}^{\mathbf{\bar{a}}}|^{2}]^{1/2}\\
&  =\sum_{k=0}^{N-1}\mathbb{E[}|\xi_{k}^{N}|]\sigma_{\mathbf{aa}^{T}}%
^{1/2}\sigma_{\mathbf{\bar{a}\bar{a}}^{T}}^{1/2}(t_{k+1}^{N}-t_{k}^{N}%
)^{3/2}\rightarrow0.
\end{align*}

\end{remark}

We now can claim our $G$--It\^{o}'s formula. Consider%
\[
X_{t}^{\nu}=X_{0}^{\nu}+\int_{0}^{t}\alpha_{s}^{\nu}ds+\int_{0}^{t}\eta
_{s}^{\nu ij}d\left \langle B^{i},B^{j}\right \rangle _{s}+\int_{0}^{t}\beta
_{s}^{\nu j}dB_{s}^{j}%
\]

\begin{proposition}
\label{Prop-Ito}Let $\alpha^{\nu}$, $\beta^{\nu j}$ and $\eta^{\nu ij}$,
$\nu=1,\cdots,n$, $i,j=1,\cdots,d$ be bounded processes of $M_{G}^{2}(0,T)$.
Then for each $t\geq0$ and in $L_{G}^{2}(\mathcal{H}_{t})$ we have%
\begin{align}
\Phi(X_{t})-\Phi(X_{s})  &  =\int_{s}^{t}\partial_{x^{\nu}}\Phi(X_{u}%
)\beta_{u}^{\nu j}dB_{u}^{j}+\int_{s}^{t}\partial_{x_{\nu}}\Phi(X_{u}%
)\alpha_{u}^{\nu}du\label{Ito-form1}\\
&  +\int_{s}^{t}[\partial_{x^{\nu}}\Phi(X_{u})\eta_{u}^{\nu ij}+\frac{1}%
{2}\partial_{x^{\mu}x^{\nu}}^{2}\Phi(X_{u})\beta_{u}^{\nu i}\beta_{u}^{\nu
j}]d\left \langle B^{i},B^{j}\right \rangle _{u}\nonumber
\end{align}

\end{proposition}

\begin{proof}
We first consider the case where $\alpha$, $\eta$ and $\beta$ are step
processes of the form%
\[
\eta_{t}(\omega)=\sum_{k=0}^{N-1}\xi_{k}(\omega)\mathbf{I}_{[t_{k},t_{k+1}%
)}(t).
\]
From the above Lemma, it is clear that (\ref{Ito-form1}) holds true. Now let
\[
X_{t}^{\nu,N}=X_{0}^{\nu}+\int_{0}^{t}\alpha_{s}^{\nu,N}ds+\int_{0}^{t}%
\eta_{s}^{\nu ij,N}d\left \langle B^{i},B^{j}\right \rangle _{s}+\int_{0}%
^{t}\beta_{s}^{\nu j,N}dB_{s}^{j}%
\]
where $\alpha^{N}$, $\eta^{N}$ and $\beta^{N}$ are uniformly bounded step
processes that converge to $\alpha$, $\eta$ and $\beta$ in $M_{G}^{2}(0,T)$ as
$N\rightarrow \infty$. From Lemma \ref{Lem-26}%
\begin{align}
\Phi(X_{t}^{N})-\Phi(X_{0})  &  =\int_{0}^{t}\partial_{x^{\nu}}\Phi(X_{u}%
^{N})\beta_{u}^{\nu j,N}dB_{u}^{j}+\int_{0}^{t}\partial_{x_{\nu}}\Phi
(X_{u}^{N})\alpha_{u}^{\nu,N}du\label{N-Ito}\\
&  +\int_{0}^{t}[\partial_{x^{\nu}}\Phi(X_{u}^{N})\eta_{u}^{\nu ij,N}+\frac
{1}{2}\partial_{x^{\mu}x^{\nu}}^{2}\Phi(X_{u}^{N})\beta_{u}^{\mu i,N}\beta
_{u}^{\nu j,N}]d\left \langle B^{i},B^{j}\right \rangle _{u}\nonumber
\end{align}
Since%
\begin{align*}
&  \mathbb{E[}|X_{t}^{N,\mu}-X_{t}^{\mu}|^{2}]\\
&  \leq C\int_{0}^{T}\{ \mathbb{E}[(\alpha_{s}^{\mu,N}-\alpha_{s}^{\mu}%
)^{2}]+\mathbb{E}[|\eta_{s}^{\mu,N}-\eta_{s}^{\mu}|^{2}]+\mathbb{E}[(\beta
_{s}^{\mu,N}-\beta_{s}^{\mu})^{2}]\}ds\\
&
\end{align*}
We then can prove that, in $M_{G}^{2}(0,T)$,
\begin{align*}
\partial_{x^{\nu}}\Phi(X_{\cdot}^{N})\eta_{\cdot}^{\nu ij,N}  &
\rightarrow \partial_{x^{\nu}}\Phi(X_{\cdot})\eta_{\cdot}^{\nu ij}\\
\partial_{x^{\mu}x^{\nu}}^{2}\Phi(X_{\cdot}^{N})\beta_{\cdot}^{\mu i,N}%
\beta_{\cdot}^{\nu j,N}  &  \rightarrow \partial_{x^{\mu}x^{\nu}}^{2}%
\Phi(X_{\cdot})\beta_{\cdot}^{\mu i}\beta_{\cdot}^{\nu j}\\
\partial_{x_{\nu}}\Phi(X_{\cdot}^{N})\alpha_{\cdot}^{\nu,N}  &  \rightarrow
\partial_{x_{\nu}}\Phi(X_{\cdot})\alpha_{\cdot}^{\nu}\\
\partial_{x^{\nu}}\Phi(X_{\cdot}^{N})\beta_{\cdot}^{\nu j,N}  &
\rightarrow \partial_{x^{\nu}}\Phi(X_{\cdot})\beta_{\cdot}^{\nu j}%
\end{align*}
We then can pass limit in both sides of (\ref{N-Ito}) and thus prove
(\ref{Ito-form1}).
\end{proof}

\section{$G$--martingales, $G$--convexity and Jensen's inequality}

\subsection{The notion of $G$--martingales}

We now give the notion of $G$--martingales:

\begin{definition}
A process $(M_{t})_{t\geq0}$ is called a $G$\textbf{--martingale} (resp.
$G$\textbf{--supermartingale}, $G$\textbf{--submartingale}) if for each $0\leq
s\leq t<\infty$, we have $M_{t}\in L_{G}^{1}(\mathcal{H}_{t})$ and
\[
\mathbb{E}[M_{t}|\mathcal{H}_{s}]=M_{s},\  \  \  \text{(resp.\ }\leq
M_{s},\  \  \geq M_{s}).
\]

\end{definition}

It is clear that, for a fixed $X\in L_{G}^{1}(\mathcal{H})$, $\mathbb{E}%
[X|\mathcal{H}_{t}]_{t\geq0}$ is a $G$--martingale. In general, how to
characterize a $G$--martingale or a $G$--supermartingale is still a very
interesting problem. But the following example gives an important characterization:

\begin{example}
Let $M_{0}\in \mathbb{R}$, $\phi=(\phi^{i})_{i=1}^{d}\in M_{G}^{2}%
(0,T;\mathbb{R}^{d})$ and $\eta=(\eta^{ij})_{i,j=1}^{d}\in M_{G}%
^{2}(0,T;\mathbb{S}_{d})$ be given and let%
\[
M_{t}=M_{0}+\int_{0}^{t}\phi_{u}^{i}dB_{s}^{j}+\int_{0}^{t}\eta_{u}%
^{ij}d\left \langle B^{i},B^{j}\right \rangle _{u}-\int_{0}^{t}2G(\eta
_{u})du,\ t\in \lbrack0,T].
\]
Then $M$ is a $G$--martingale on $[0,T]$. To prove this it suffices to prove
the case $\eta \in M_{G}^{2,0}(0,T;\mathbb{S}_{d})$, i.e.,
\[
\eta_{t}=\sum_{k=0}^{N-1}\eta_{t_{k}}I_{[t_{k}.t_{k+1})}(t).
\]
We have, for $s\in \lbrack t_{N-1},t_{N}]$,%
\begin{align*}
\mathbb{E}[M_{t}|\mathcal{H}_{s}]  &  =M_{s}+\mathbb{E}[\eta_{t_{N-1}}%
^{ij}(\left \langle B^{i},B^{j}\right \rangle _{t}-\left \langle B^{i}%
,B^{j}\right \rangle _{s})-2G(\eta_{t_{N-1}})(t-s)|\mathcal{H}_{s}]\\
&  =M_{s}+\mathbb{E}[\eta_{t_{N-1}}^{ij}(B_{t}^{i}-B_{s}^{i})(B_{t}^{j}%
-B_{s}^{j})|\mathcal{H}_{s}]-2G(\eta_{t_{N-1}})(t-s)\\
&  =M_{s}.
\end{align*}
In the last step, we apply the relation (\ref{eq-GMB-14a}). We then can repeat
this procedure, step by step backwardly, to prove the case $s\in
\lbrack0,t_{N-1}]$. \ 
\end{example}

\begin{remark}
It is worth to mention that for a $G$--martingale, in general, $-M$ is not a
$G$--martingale. But in the above example, when $\eta \equiv0$, then $-M$ is
still a $G$--martingale. This makes an essential difference of the $dB$ part
and the $d\left \langle B\right \rangle $ part of a $G$--martingale.
\end{remark}

\subsection{$G$--convexity and Jensen's inequality for $G$--expectation}

A very interesting question is whether the well--known Jensen's inequality
still holds for $G$--expectation. In the framework of $g$--expectation, this
problem was investigated in \cite{BCHMP1} in which a counterexample is given
to indicate that, even for a linear function which is obviously convex,
Jensen's inequality for $g$-expectation generally does not hold. Stimulated by
this example, \cite{JC1} proved that Jensen's inequality holds for any convex
function under a $g$--expectation if and only if the corresponding generating
function $g=g(t,z)$ is super-homogeneous in $z$. Here we will discuss this
problem in a quite different point of view. We will define a new notion of convexity:

\begin{definition}
A $C^{2}$-function $h:\mathbb{R\longmapsto R}$ is called $G$\textbf{--convex}
if the following condition holds for each $(y,z,A)\in \mathbb{R}\times
\mathbb{R}^{d}\times \mathbb{S}_{d}$:
\begin{equation}
G(h^{\prime}(y)A+h^{\prime \prime}(y)zz^{T})-h^{\prime}(y)G(A)\geq
0,\  \label{G-conv}%
\end{equation}
where $h^{\prime}$ and $h^{\prime \prime}$ denote the first and the second
derivatives of $h$.
\end{definition}

It is clear that in the special situation where $G(D^{2}u)=\frac{1}{2}\Delta
u$, the corresponding $G$-convex function becomes a classical convex function.

\begin{lemma}
The following two conditions are equivalent:\newline \textbf{(i)} the function
$h$ is $G$--convex.\newline \textbf{(ii)} The following Jensen inequality
holds: for each $T\geq0$,
\begin{equation}
\mathbb{E}[h(\phi(B_{T}))]\geq h(\mathbb{E}[\phi(B_{T})]), \label{gg-Jensen}%
\end{equation}
for each $C^{2}$--function $\phi$ such that $h(\phi(B_{T}))$, $\phi(B_{T})\in
L_{G}^{1}(\mathcal{H}_{T})$.
\end{lemma}

\begin{proof}
(i) $\Longrightarrow$(ii) By the definition $u(t,x):=P_{t}^{G}[\phi
](x)=\mathbb{E}[\phi(x+B_{t})]$ solves the nonlinear heat equation
(\ref{eq-heat}). Here we only consider the case where $u$ is a $C^{1,2}%
$-function. Otherwise we can use the language of viscosity solution as we did
in the proof of Lemma \ref{Scaling}. By simple calculation, we have
\[
\partial_{t}h(u(t,x))=h^{\prime}(u)\partial_{t}u=h^{\prime}(u(t,x))G(D^{2}%
u(t,x)),
\]
or%
\[
\partial_{t}h(u(t,x))-G(D^{2}h(u(t,x)))-f(t,x)=0,\ h(u(0,x))=h(\phi(x)),
\]
where we denote%
\[
f(t,x)=h^{\prime}(u(t,x))G(D^{2}u(t,x))-G(D^{2}h(u(t,x))).
\]
Since $h$ is $G$--convex, it follows that $f\leq0$ and thus $h(u)$ is a
$G$-subsolution. It follows from the maximum principle that $h(P_{t}^{G}%
(\phi)(x))\leq P_{t}^{G}(h(\phi))(x)$. In particular (\ref{gg-Jensen}) holds.
Thus we have (ii).\newline(ii) $\Longrightarrow$(i): For a fixed
$(y,z,A)\in \mathbb{R\times R}^{d}\times \mathbb{S}_{d}$, we set $\phi
(x):=y+\left \langle x,z\right \rangle +\frac{1}{2}\left \langle
Ax,x\right \rangle $. By the definition of $P_{t}^{G}$ we have $\partial
_{t}(P_{t}^{G}(\phi)(x))|_{t=0}=G(D^{2}\phi)(x)$. By (ii) we have%
\[
h(P_{t}^{G}(\phi)(x))\leq P_{t}^{G}(h(\phi))(x).
\]
Thus, for $t>0$,%
\[
\frac{1}{t}[h(P_{t}^{G}(\phi)(x))-h(\phi(x))]\leq \frac{1}{t}[P_{t}^{G}%
(h(\phi))(x)-h(\phi(x))]
\]
We then let $t$ tend to $0$:%
\[
h^{\prime}(\phi(x))G(D^{2}\phi(x))\leq G(D_{xx}^{2}h(\phi(x))).
\]
Since $D_{x}\phi(x)=z+Ax$ and $D_{xx}^{2}\phi(x)=A$. We then set $x=0$ and
obtain (\ref{G-conv}).
\end{proof}

\begin{proposition}
The following two conditions are equivalent:\newline \textbf{(i)} the function
$h$ is $G$--convex.\newline \textbf{(ii)} The following Jensen inequality
holds:
\begin{equation}
\mathbb{E}[h(X)|\mathcal{H}_{t}]\geq h(\mathbb{E}[X|\mathcal{H}_{t}%
]),\  \ t\geq0, \label{JensenX}%
\end{equation}
for each $X\in L_{G}^{1}(\mathcal{H})$ such that $h(X)\in L_{G}^{1}%
(\mathcal{H})$.
\end{proposition}

\begin{proof}
The part (ii) $\Longrightarrow$(i) is already provided by the above Lemma. We
can also apply this lemma to prove (\ref{JensenX}) for the case $X\in
L_{ip}^{0}(\mathcal{H})$ of the form $X=\phi(B_{t_{1}},\cdots,B_{t_{m}%
}-B_{t_{m-1}})$ by using the procedure of the definition of $\mathbb{E}%
[\cdot]$ and $\mathbb{E}[\cdot|\mathcal{H}_{t}]$ given in Definition
\ref{Def-3}. We then can extend this Jensen's inequality, under the norm
$\left \Vert \cdot \right \Vert =\mathbb{E}[|\cdot|]$, to the general situation.
\end{proof}

\begin{remark}
The above notion of $G$--convexity can be also applied to the case where the
nonlinear heat equation (\ref{eq-heat}) has a more general form:
\begin{equation}
\frac{\partial u}{\partial t}-G(u,\nabla u,D^{2}u)=0,\  \ u(0,x)=\psi(x)
\label{PG-psi}%
\end{equation}
(see Examples 4.3, 4.4 and 4.5 in \cite{Peng2005}). In this case a $C^{2}%
$-function $h:\mathbb{R\longmapsto R}$ is called to be $G$--\textbf{convex} if
the following condition holds for each $(y,z,A)\in \mathbb{R}\times
\mathbb{R}^{d}\times \mathbb{S}_{d}$:
\[
G(y,h^{\prime}(y)z,h^{\prime}(y)A+h^{\prime \prime}(y)zz^{T})-h^{\prime
}(y)G(y,z,A)\geq0.
\]
We don't need the subadditivity and/or positive homogeneity of $G(y,z,A)$. A
particularly interesting situation is the case of $g$--expectation for a given
generating function $g=g(y,z)$, $(y,z)\in \mathbb{R}\times \mathbb{R}^{d}$, in
this case we have the following $g$--convexity:
\begin{equation}
\frac{1}{2}h^{\prime \prime}(y)|z|^{2}+g(h(y),h^{\prime}(y)z)-h^{\prime
}(y)g(y,z)\geq0. \label{g-convex}%
\end{equation}
We will discuss this situation elsewhere.
\end{remark}

\begin{example}
Let $h$ be a $G$--convex function and let $X\in L_{G}^{1}(\mathcal{H})$ be
such that $h(X)\in L_{G}^{1}(\mathcal{H})$, then $Y_{t}=h(\mathbb{E}%
[X|\mathcal{H}_{t}])$, $t\geq0$, is a $G$--submartingale: for each $s\leq t$,
\[
\mathbb{E}[Y_{t}|\mathcal{H}_{s}]=\mathbb{E}[h(\mathbb{E}[X|\mathcal{F}%
_{t}])|\mathcal{F}_{s}]\geq h(\mathbb{E}[X|\mathcal{F}_{s}])=Y_{s}\text{.}%
\]

\end{example}

\section{Stochastic differential equations}

We consider the following SDE driven by $G$-Brownian motion.%
\begin{equation}
X_{t}=X_{0}+\int_{0}^{t}b(X_{s})ds+\int_{0}^{t}h_{ij}(X_{s})d\left \langle
B^{i},B^{j}\right \rangle _{s}+\int_{0}^{t}\sigma_{j}(X_{s})dB_{s}^{j}%
,\ t\in \lbrack0,T]. \label{SDE}%
\end{equation}
where the initial condition $X_{0}\in \mathbb{R}^{n}$ is given and
\[
b,h_{ij},\sigma_{j}:\mathbb{R}^{n}\mapsto \mathbb{R}^{n}%
\]
are given Lipschitz functions, i.e., $|\phi(x)-\phi(x^{\prime})|\leq
K|x-x^{\prime}|$, for each $x$, $x^{\prime}\in \mathbb{R}^{n}$, $\phi=b$,
$\eta_{ij}$ and $\sigma_{j}$. Here the horizon $[0,T]$ can be arbitrarily
large. The solution is a process $X\in M_{G}^{2}(0,T;\mathbb{R}^{n})$
satisfying the above SDE. We first introduce the following mapping on a fixed
interval $[0,T]$:%
\[
\Lambda_{\cdot}(Y):=:Y\in M_{G}^{2}(0,T;\mathbb{R}^{n})\longmapsto M_{G}%
^{2}(0,T;\mathbb{R}^{n})\  \
\]
by setting $\Lambda_{t}=X_{t}$, $t\in \lbrack0,T]$, with
\[
\Lambda_{t}(Y)=X_{0}+X_{0}+\int_{0}^{t}b(Y_{s})ds+\int_{0}^{t}h_{ij}%
(Y_{s})d\left \langle B^{i},B^{j}\right \rangle _{s}+\int_{0}^{t}\sigma
_{j}(Y_{s})dB_{s}^{j}.
\]

We immediately have

\begin{lemma}
For each $Y,Y^{\prime}\in M_{G}^{2}(0,T;\mathbb{R}^{n})$, we have the
following estimate:%
\[
\mathbb{E}[|\Lambda_{t}(Y)-\Lambda_{t}(Y^{\prime})|^{2}]\leq C\int_{0}%
^{t}\mathbb{E}[|Y_{s}-Y_{s}^{\prime}|^{2}]ds,\ t\in \lbrack0,T],
\]
where the constant $C$ depends only on $K$, $\Gamma$ and the dimension $n$.
\end{lemma}

\begin{proof}
This is a direct consequence of the inequalities (\ref{Bohner}), (\ref{e2})
and (\ref{qua-ine}).
\end{proof}

We now prove that SDE (\ref{SDE}) has a unique solution. By multiplying
$e^{-2Ct}$ on both sides of the above inequality and then integrate them on
$[0,T]$. It follows that%
\begin{align*}
\int_{0}^{T}\mathbb{E}[|\Lambda_{t}(Y)-\Lambda_{t}(Y^{\prime})|^{2}%
]e^{-2Ct}dt  &  \leq C\int_{0}^{T}e^{-2Ct}\int_{0}^{t}\mathbb{E}[|Y_{s}%
-Y_{s}^{\prime}|^{2}]dsdt\\
&  =C\int_{0}^{T}\int_{s}^{T}e^{-2Ct}dt\mathbb{E}[|Y_{s}-Y_{s}^{\prime}%
|^{2}]ds\\
&  =(2C)^{-1}C\int_{0}^{T}(e^{-2Cs}-e^{-2CT})\mathbb{E}[|Y_{s}-Y_{s}^{\prime
}|^{2}]ds.
\end{align*}
We then have
\[
\int_{0}^{T}\mathbb{E}[|\Lambda_{t}(Y)-\Lambda_{t}(Y^{\prime})|^{2}%
]e^{-2Ct}dt\leq \frac{1}{2}\int_{0}^{T}\mathbb{E}[|Y_{t}-Y_{t}^{\prime}%
|^{2}]e^{-2Ct}dt.
\]
We observe that the following two norms are equivalent in $M_{G}%
^{2}(0,T;\mathbb{R}^{n})$
\[
\int_{0}^{T}\mathbb{E}[|Y_{t}|^{2}]dt\thicksim \int_{0}^{T}\mathbb{E}%
[|Y_{t}|^{2}]e^{-2Ct}dt.
\]
From this estimate we can obtain that $\Lambda(Y)$ is a contract mapping.
Consequently, we have

\begin{theorem}
There exists a unique solution of $X\in M_{G}^{2}(0,T;\mathbb{R}^{n})$ of the
stochastic differential equation (\ref{SDE}).
\end{theorem}

\section{Appendix: Some inequalities in $L_{G}^{p}(\mathcal{H})$}

For $r>0$, $1<p,q<\infty$, such that $\frac{1}{p}+\frac{1}{q}=1$, we have
\begin{align}
|a+b|^{r}  &  \leq \max \{1,2^{r-1}\}(|a|^{r}+|b|^{r}),\  \  \forall
a,b\in \mathbb{R}\label{ee4.3}\\
|ab|  &  \leq \frac{|a|^{p}}{p}+\frac{|b|^{q}}{q}. \label{ee4.4}%
\end{align}

\begin{proposition}%
\begin{align}
\mathbb{E}[|X+Y|^{r}]  &  \leq C_{r}(\mathbb{E}[|X|^{r}]+\mathbb{E[}%
|Y|^{r}]),\label{ee4.5}\\
\mathbb{E}[|XY|]  &  \leq \mathbb{E}[|X|^{p}]^{1/p}\cdot \mathbb{E}%
[|Y|^{q}]^{1/q},\label{ee4.6}\\
\mathbb{E}[|X+Y|^{p}]^{1/p}  &  \leq \mathbb{E}[|X|^{p}]^{1/p}+\mathbb{E}%
[|Y|^{p}]^{1/p} \label{ee4.7}%
\end{align}
In particular, for $1\leq p<p^{\prime}$, we have $\mathbb{E}[|X|^{p}%
]^{1/p}\leq \mathbb{E}[|X|^{p^{\prime}}]^{1/p^{\prime}}.$
\end{proposition}

\begin{proof}
(\ref{ee4.5}) follows from (\ref{ee4.3}). We set
\[
\xi=\frac{X}{\mathbb{E}[|X|^{p}]^{1/p}},\  \  \eta=\frac{Y}{\mathbb{E}%
[|Y|^{q}]^{1/q}}.
\]
By (\ref{ee4.4}) we have%
\begin{align*}
\mathbb{E}[|\xi \eta|]  &  \leq \mathbb{E}[\frac{|\xi|^{p}}{p}+\frac{|\eta|^{q}%
}{q}]\leq \mathbb{E}[\frac{|\xi|^{p}}{p}]+\mathbb{E}[\frac{|\eta|^{q}}{q}]\\
&  =\frac{1}{p}+\frac{1}{q}=1.
\end{align*}
Thus (\ref{ee4.6}) follows.
\begin{align*}
\mathbb{E}[|X+Y|^{p}]  &  =\mathbb{E}[|X+Y|\cdot|X+Y|^{p-1}]\\
&  \leq \mathbb{E}[|X|\cdot|X+Y|^{p-1}]+\mathbb{E}[|Y|\cdot|X+Y|^{p-1}]\\
&  \leq \mathbb{E}[|X|^{p}]^{1/p}\cdot \mathbb{E[}|X+Y|^{(p-1)q}]^{1/q}\\
&  +\mathbb{E}[|Y|^{p}]^{1/p}\cdot \mathbb{E[}|X+Y|^{(p-1)q}]^{1/q}%
\end{align*}
We observe that $(p-1)q=p$. Thus we have (\ref{ee4.7}).
\end{proof}

\end{document}